\documentclass[12pt]{report} 
\usepackage{latexsym}
\newcommand{\ds}{\displaystyle}

\usepackage{amsfonts}

\usepackage{epsfig}

\newcommand{\reg}{\mbox{\scriptsize{reg}}}

\newcommand{\codim}{\mbox{codim}}
\newcommand{\trdeg}{\mbox{trdeg}}
\newcommand{\sing}{\mbox{\scriptsize{sing}}}

\newcommand{\rank}{\mbox{rank}}

\newcommand{{\CB}}{\cal B}

\newcommand{{\CC}}{\cal C}
\def\subs{\mathop{\subset}}

 1
 1

\def\picture#1 by #2 (#3){
\vbox to #2{
\hrule width #1 height 0pt depth 0pt
\vfill
\special{picture #3}}}

\begin{document}
\begin{center}
{\bf \Large The maximal abelian dimension of a Lie algebra, Rentschler's property and Milovanov's conjecture}\\
\ \\
{\bf \large Alfons I. Ooms}\\
\ \\
{\it Mathematics Department, Hasselt University, Agoralaan, Campus Diepenbeek, 3590 Diepenbeek, Belgium\\
E-mail address: alfons.ooms@uhasselt.be}\\
To the memory of Maryam Mirzakhany (1977-2017)
\end{center}

{\bf Key words:} maximal abelian dimension, Rentschler's property, complete Poisson commutative subalgebras, filiform Lie algebras, Milovanov's conjecture\\
MSC: 17B30, 17B63.\\
\ \\
{\bf Abstract.}\\
A finite dimensional Lie algebra $L$ with magic number $c(L)$ is said to satisfy Rentschler's property if it admits an abelian Lie subalgebra $H$ of dimension at least $c(L)-1$.  
We study the occurrence of this new property in various Lie algebras, such as nonsolvable, solvable, nilpotent, metabelian and filiform Lie algebras.  Under some mild condition $H$ gives rise to a complete Poisson commutative
subalgebra of the symmetric algebra $S(L)$. Using this, we show that Milovanov's conjecture holds for the filiform Lie algebras of type $L_n$, $Q_n$, $R_n$, $W_n$ and also for all filiform Lie algebras of dimension at most eight. 
For the latter the Poisson center of these Lie algebras is determined.\\
\ \\
{\bf \large 1. Introduction}\\
Let $k$ be an algebraically closed field of characteristic zero and let $L$ be a Lie algebra over $k$ with basis $x_1, \ldots, x_n$. Let $S(L) = k[x_1,\ldots, x_n]$ be its symmetric algebra with quotient
field $R(L)$. For each $\xi \in L^\ast$ we consider its stabilizer
\begin{eqnarray*}
L(\xi) = \{x \in L\mid \xi ([x,y]) = 0\ \mbox{for all}\ y \in L\}
\end{eqnarray*}
The minimal value of $\dim L(\xi)$ is called the index of $L$ and is denoted by $i(L)$ [D, 1.11.6; TY, 19.7.3].  
An element $\xi \in L^\ast$ is called regular if $\dim L(\xi) = i(L)$.  The set $L^\ast_{\reg}$ of all regular elements of $L^\ast$ is an open dense subset of $L^\ast$.\\
We put $L^\ast_{\sing} = L^\ast \backslash L^\ast_{\reg}$.  Clearly, $\codim\ L^\ast_{\sing} \geq 1$.  Following [JS] we call $L$ singular if equality holds and nonsingular otherwise.  
For instance, any semi-simple Lie algebra $L$ is nonsingular since $\codim\ L^\ast_{\sing} = 3$. We recall from [D, 1.14.13] that
\begin{eqnarray}
i(L) = \dim L - \rank_{R(L)} ([x_i,x_j])\label{1}
\end{eqnarray}
In particular, $\dim L - i(L)$ is an even number.\\
The integer $c(L) = (\dim L + i(L))/2$ is called the magic number of $L$ and in this paper it will certainly live up to its name. Also, the Frobenius semiradical $F(L)$ of $L$ and quasi quadratic Lie algebras 
will play an important role (see Section 2 for definitions, properties and examples).\\
Next, we denote by $p_L \in S(L)$ the fundamental semi-invariant of $L$, i.e. the greatest common divisor of the Pfaffians of the principal $t \times t$ minors of the structure matrix 
$B = ([x_i,x_j]) \in M_n(R(L))$ where $t = \rank\ B$ 
[O6, 2.6].  It is well-known that $L$ is singular if and only is $p_L \notin k$ [OV, p.307].\\
\ \\
{\bf 1.1 The Poisson algebra $\boldmath{S(L)}$ and its center}\\
The symmetric algebra $S(L)$, which we identify with $k[x_1,\ldots, x_n]$, has a natural Poisson algebra structure, the Poisson bracket of $f, g \in S(L)$ given by:
\begin{eqnarray}
\{f,g\} = \sum\limits_{i=1}^n \sum\limits_{j=1}^n [x_i,x_j] \ds\frac{\partial f}{\partial x_i} \ds\frac{\partial g}{\partial x_j}\label{2}
\end{eqnarray}
The Poisson center $Y(L)$ of $S(L)$ coincides with the subalgebra $S(L)^L$ of invariants of $S(L)$.  We say that $L$ is coregular if $Y(L)$ is a polynomial algebra over $k$.\\
Let $A$ be a Poisson commutative subalgebra of $S(L)$ (i.e. $\{f,g\} = 0$ for all $f, g \in A$).  Then it is well-known that $\trdeg_k(A) \leq c(L)$.  
$A$ is called complete if equality holds and strongly complete if it is also a maximal Poisson commutative subalgebra of $S(L)$.  According to Sadetov there always exists a complete Poisson commutative 
subalgebra of $S(L)$ [Sa].\\
\ \\
In 1999 Milovanov raised the following interesting open question, which has been verified for some low-dimensional Lie algebras [Ko1], in particular for all complex nilpotent Lie algebras of dimension at most seven
[O5], see also [Ko2]. In this paper the conjecture is shown to hold for the filiform Lie algebras of type $L_n$, $Q_n$, $R_n$, $W_n$ (Theorem 8.6) as well as for all filiform Lie algebras of dimension
$\leq 8$ (Corollary 9.3) It is also verified for various solvable Lie algebras we will come across.\\
\newpage
{\bf The Milovanov conjecture.}\\
For any solvable Lie algebra $L$ there exists a complete Poisson commutative subalgebra $M$ of $S(L)$ generated by elements of degree at most two.\\
\ \\
{\bf 1.2 The maximal abelian dimension $\boldmath{\alpha(L)}$}\\
This is by definition the maximum dimension of all abelian subalgebras of $L$.  Note that the dimension of a maximal abelian subalgebra of $L$ may be strictly smaller than $\alpha(L)$ [O2, Example 1], [BC, Example 2.2].  Furthermore, if $L$ is solvable then there is always an abelian
ideal of $L$ with dimension $\alpha(L)$ [BC, Proposition 2.6].\\
\ \\
$\alpha(L)$ is an interesting invariant for $L$.  Already in 1905 Schur proved that 
$$\alpha(gl(n,\Bbb C)) = [n^2/4]+1$$ 
where [\ ] is the greatest integer function [Sc, p.67].  
In 1944 Jacobson generalized this by replacing $\Bbb C$ by any field [J, p.434].  In 1998 Mirzakhani, gave a one-page proof, by induction on $n$, of  Schur's result [Mir].  We will also provide a short proof by making use of the magic number (Corollary 6.10).  
See also [Co, G, K, SVVZ, WL].  In 1945 Malcev completely determined the abelian subalgebras of semi-simple Lie algebras
[Ma]. 
\begin{center}
{\bf Table 1. The invariant $\boldmath{\alpha}$ for a simple Lie algebra $\boldmath{L}$}\\
\ \\
\begin{tabular}{ccc}
\hline
$L$ &$\dim L$ &$\alpha(L)$\\
\hline
$A_n, n \geq 1$ &$n(n+2)$ &$\lceil \left(\frac{n+1}{2}\right)^2\rceil$\\
\hline
$B_3$ &21 &5\\
\hline
$B_n, n \geq 4$ &$n(2n+1)$ &$\frac{n(n-1)}{2} + 1$\\
\hline
$C_n, n \geq 2$ &$n(2n+1)$ &$\frac{n(n+1)}{2}$\\
\hline
$D_n, n \geq 4$ &$n(2n-1)$ &$\frac{n(n-1)}{2}$\\
\hline
$G_2$ &14 &3\\
\hline
$F_4$ &52 &9\\
\hline
$E_6$ &78 &16\\
\hline
$E_7$ &133 &27\\
\hline
$E_8$ &248 &36\\
\hline
\end{tabular}
\end{center}
Later on, Kostant [K], Suter [Su] and Panyushev [P] studied the abelian ideals of a Borel subalgebra of $L$.\\
Over the years the invariant $\alpha(L)$ and its general properties have drawn much attention for various reasons.  Several bounds on $\alpha(L)$ have been found depending on the nature of $L$.  
For instance, $(\sqrt{8n+9} - 3)/2 \leq \alpha(L)$ if $L$ is solvable and $(\sqrt{8n+1} - 1)/2 \leq \alpha(L)$ if $L$ is nilpotent, where $n = \dim L$ [BC, C, CT, Mil1, Mil2, Mil3, Mo, R, St].\\
A useful algorithm has already been established in order to compute $\alpha(L)$, in particular for all indecomposable solvable Lie algebras of dimension at most 6 and for all indecomposable nilpotent Lie algebras of dimension at most 7.  
The same technique was also employed to determine $\alpha$ for the following Lie algebras: the $n$-dimensional Heisenberg Lie algebra $H_n$ and $T_n$ (respectively $N_n$) 
the Lie algebra of all $n \times n$ upper (resp. strictly upper) triangular matrices [BNT1, BNT2, C, CNT1, CNT2, CNT3, T].  However shorter proofs can be given (Proposition 6.7, Theorem 6.9 and Corollary 6.10) by making use of the 
following simple observation, which seems to have been ignored in the literature.\\
For any Lie algebra $L$ the following holds:
\begin{eqnarray}
i(L) \leq \alpha(L) \leq c(L)\label{3}
\end{eqnarray}
To show the first inequality, just take a regular linear functional $\xi \in L^\ast$.  Its stabilizer $L(\xi)$ is an abelian subalgebra of $L$ [D1.11.7].  Hence, $i(L) = \dim L(\xi) \leq \alpha (L)$.  
See [O2, Theorem 14] for the second inequality.  
We will encounter a few examples for which $i(L) = \alpha(L)$, for instance $\mathfrak{g}_{5,4}$, $sl(2,k)$, $sl(3,k)$, Examples 5.4 and 6.4.  On the other hand, the equality $\alpha(L) = c(L)$ happens more frequently.
In that case there is a commutative subalgebra $P$ of $L$ such that $\dim P = c(L)$, in other words $P$ is a commutative polarization (CP) with respect to any regular $\xi \in L^\ast$.  
These CP's play an important role in the study of the irreducible representations of the enveloping algebra $U(L)$ and their kernels, the primitive ideals of $U(L)$ [EO].  Moreover, $S(P)$ is a strongly complete Poisson
commutative subalgebra of $S(L)$, in particular it contains the Poisson center $Y(L)$ of $S(L)$.  Clearly $S(P)$ is a polynomial algebra over $k$, generated by elements of degree one, so in this situation Milovanov's conjecture is trivially
satisfied.  Also $R(P)$ is a maximal, Poisson commutative subfield of $R(L)$ [O2].  CP's often appear in low-dimensional solvable Lie algebras and even more so in the nilpotent case.  For instance in the list of the 159 indecomposable nilpotent Lie algebras up to dimension 7
(a family is counted as one member) only 28 do not possess a CP [O4, O5].  Rudolf Rentschler discovered that each Lie algebra $L$ of these 28 exceptions still admits an 
abelian subalgebra $H$ of $L$ of dimension $c(L)-1$.  Moreover if $Y(L) \not\subset S(H)$ 
then the Poisson commutative subalgebra $M = S(H) Y(L)$ is complete in $S(L)$.  This led us to introduce the following definition.\\
\ \\
{\bf 1.3 Rentschler's property}\\
We say that a Lie algebra $L$ satisfies Rentschler's property ($R$-property for short) if $L$ either admits a CP or an abelian subalgebra $H$ of dimension $c(L)-1$.  In other words $\alpha (L) \geq c(L)-1$.\\
\ \\
{\bf Remark.} Let $H$ be an abelian subalgebra of $L$ of dimension $c(L)-1$.  Then $M = S(H) Y(L)$ is a complete Poisson commutative subalgebra of $S(L)$ if and only if $Y(L) \not\subset S(H)$ (however $M$ is not necessarily strongly complete.  See 2.2, Example 3).  
Suppose in addition that $L$ is coregular without proper semi-invariants (or more generally satisfying the Joseph-Shafrir conditions [O6, Definition 17]) and that
\begin{eqnarray*}
3i(L) + 2 \deg p_L = \dim L + 2 \dim Z(L)
\end{eqnarray*}
Then $Y(L)$ is generated by elements of degree at most 2 [O6, Corollary 19] and hence the same holds for $M$.\\
\ \\
An $n$-dimensional Lie algebra $L$ satisfies the $R$-property if
\begin{itemize}
\item[1.] $L$ is indecomposable and nonsolvable with $n \leq 7$.  There are only 3 exceptions if $n=8$ (Theorem 4.1).
\item[2.] (i) $L$ is solvable with $n \leq 6$ (Theorem 5.2 (1)).\\
(ii) $n=7$, $L$ is solvable but not quasi quadratic (Theorem 5.2(2)).  The latter condition cannot be removed (Example 5.4).
\item[3.] (i) $L$ is nilpotent with $n \leq 7$ (Theorem 6.2(1)).\\
(ii) $n=8$, $L$ is nilpotent but§ not quasi quadratic (Theorem 6.2(2)). The latter condition cannot be removed (Example 6.4).
\item[4.] $L$ is metabelian with $\dim[L,L] \geq 2$ and $n \leq 9$ (Theorem 7.4).
\item[5.] $L$ is filiform with $n \leq 11$ (Theorem 8.7) as well as those of type $L_n$, $Q_n$, $R_n$ and $W_n$ (Theorem 8.6).
\end{itemize}

The existence of CP's is not limited to low-dimensional Lie algebras [EO], so the same holds a fortiori for the $R$-property. (Lemma 3.2, Proposition 4.2, Example 5.6, Theorem 6.6, Proposition 6.7, Theorem 6.9, Corollary 7.2, Proposition 8.5
and Theorem 8.6).\\
\ \\
In section 9 we determine, case by case, the Poisson center (among other things) of all complex filiform Lie algebras $\mathfrak{g}$ of dimension at most eight, since it is a useful tool for the construction of complete Poisson
commutative subalgebras of $S(\mathfrak{g})$ (see 9.2).  This allows us to verify Milovanov's conjecture in this situation (Corollary 9.3).  We use MAPLE for the less trivial calculations.\\
\ \\
{\bf \large 2. The Frobenius semiradical $\boldmath{F(L)}$ of $\boldmath{L}$ and quasi quadratic Lie algebras [O3]}\\
Let $L$ be a Lie algebra. Its Frobenius semi-radical $F(L)$ is defined as follows:
$$F(L) = \sum\limits_{\xi \in L_{\reg}^\ast} L( \xi)$$
$F(L)$ is a characteristic ideal of $L$ containing the center $Z(L)$ of $L$ and $F(F(L)) = F(L)$.  Moreover,
$$F(L) = 0 \ \Leftrightarrow\ L \ \mbox{is Frobenius (i.e.}\ i(L) = 0\ [O1])$$
$F(L)$ can also be characterized as follows:\\
The Poisson center $Y(L)$ is contained in $S(F(L))$ and $F(L)$ is the smallest subalgebra of $L$ with this property in case $L$ is an algebraic Lie algebra without proper semi-invariants.\\
If $L$ admits a CP $P$ then $F(L) \subset P$ and hence $F(L)$ is commutative in this situation.\\
\ \\
We call $L$ quasi quadratic if $F(L) = L$. Such Lie algebras have no proper semi-invariants.  In particular they are unimodular, $\trdeg_kY(L) = i(L)$ and $R(L)^L = Q(Y(L))$ (i.e. the quotient field of $Y(L))$.  They form a large class containing
all quadratic Lie algebras (i.e. those equipped with a nondegenerate, invariant symmetric bilinear form), in particular all semi-simple and abelian Lie algebras.  However in the nilpotent case they
are rather rare (for instance there are only 5 indecomposable ones in dimension 7 [O5, p.111]).\\
\ \\
Finally, assume that $L$ is quasi quadratic.
\begin{itemize}
\item [(i)] If $L$ is solvable and nonzero, then $Z(L) \neq 0$.
\item[(ii)] If $L$ is nilpotent with $\dim L \geq 2$ then $\dim Z(L) \geq 2$.
\end{itemize}

{\bf 2.1 Remark}\\
A number of years ago Mustaphe Rais kindly sent us an unpublished manuscript by Andr\'e Cerezo [Cer], in which the soul (resp. the rational soul) of a Lie algebra $L$ is introduced and studied.
This is the smallest Lie subalgebra of $L$ whose enveloping algebra (resp. enveloping quotient division ring) contains the center $Z(U(L))$ (resp. $Z(D(L))$).  
It is not hard to see that the rational soul of $L$ coincides with $F(L)$ in case $L$ is algebraic.\\
\ \\
{\bf 2.2 Examples}\\
Each of the following Lie algebras $L$ is quasi quadratic.  Only the first two are quadratic.  Since $F(L) = L$ is not commutative, $L$ does not admit a CP.  
Hence $\alpha(L) \leq c(L) -1$ by formula (\ref{2}).  By exhibiting an abelian subalgebra $H$ of $L$ of dimension $c(L)-1$, we may conclude that $L$ satisfies the $R$-property and also that 
$\alpha(L) = c(L) - 1$.  Below we will use the following abbreviations: $i= i(L)$, $c=c(L)$, $\alpha = \alpha(L)$, $Y = Y(L)$, $p = p_L$ and $M$ will be a complete, Poisson commutative subalgbra of $S(L)$.
\begin{itemize}
\item[(1)] $D_4$, the Diamond Lie algebra with basis $t, x, y, z$ and nonzero brackets:
$$[t,x] = x, [t,y] = - y, [x,y] = z$$
$D_4$ is solvable and quadratic.  It is the smallest solvable, nonabelian quasi quadratic Lie algebra. $D_4$ will be generalized in Example 5.6.
\item[(2)] $L = \mathfrak{g}_{5.4}$ with basis $x_1, x_2, x_3, x_4, x_5$ and nonzero brackets
$$[x_1, x_2] = x_3, [x_1, x_3] = x_4, [x_2, x_3] = x_5$$
$L$ is nilpotent and quadratic.\\
$i=3$, $c=4$, $p=1$, $H = \langle x_3, x_4, x_5\rangle$ is abelian, $\alpha = 3 = i = c-1$.  $ Y = k[x_4, x_5, x_3^2 + 2x_1x_5 - 2x_2x_4]$, $M = S(H)Y = k[x_3, x_4, x_5, x_1x_5 - x_2x_4]$ 
is strongly complete.
\item[(3)] Let $L$ be the 6-dimensional Lie algebra with basis $x_1,x_2,\ldots, x_6$ and nonzero brackets:\\
$[x_1,x_2] = -x_2$, $[x_1, x_3] = 2x_3$, $[x_1, x_4]  = -2x_4$, $[x_1,x_5] = x_5$,\\
$[x_2, x_3] = x_5$, $[x_2, x_5] = x_6$, $[x_3, x_4] = x_6$.\\
$i=2$, $c=4$, $p=1$.  $L$ is solvable and $H = \langle x_4, x_5, x_6\rangle$ is abelian, $\alpha = 3 = c-1$.  $Y = k [x_6, f]$ where
$$f = x_4x_5^2 + x_2x_5x_6 - 2x_3x_4x_6 - x_1x_6^2$$
$$M = S(H)Y = k [x_4, x_5, x_6, (x_2x_5 - 2x_3x_4 - x_1x_6)x_6]$$
which is complete but not strongly complete while the following
$$M_1 = k [x_4, x_5,x_6, x_2x_5 - 2x_3x_4 - x_1x_6]$$
is strongly complete and the Milovanov conjecture holds.
\item[(4)] Let $L$ be the 7-dimensional Lie algebra with basis $x_1,x_2,\ldots,x_7$ and nonzero brackets:\\
$[x_1,x_2] = x_2$, $[x_1,x_3] = -x_3$, $[x_1,x_5] = -x_5$, $[x_1,x_6] = x_6$, \\
$[x_2, x_4] = x_6$, $[x_2,x_5] = x_7$, $[x_3, x_4] = -x_5$, $[x_3, x_6] = -x_7$\\
$i=3$, $c=5$, $p=1$\\
$L$ is solvable and $H = \langle x_4, x_5, x_6, x_7\rangle$ is abelian, $\alpha = 4 = c - 1$.\\
$Y = k [x_7, h,g]$, $f = x_5x_6 - x_4x_7$, $g = x_1x_7 + x_2x_5 + x_3x_6$
$$M = S(H)Y = k[x_4,x_5,x_6,x_7,g]$$
\item[(5)] $L = L_{6,3}$ (see [O3; O6, Example 54]) with basis $h,x,y,e_0, e_1, e_2$ and nonzero brackets:\\
$[h,x] = 2x$, $[h,y] = -2y$, $[x,y] = h$, $[h,e_0] = e_0$,\\
$[h,e_1] = -e_1$, $[x,e_1] = e_0$, $[y,e_0] = e_1$, $[e_0, e_1] = e_2$.\\
$i=2$, $c=4$, $p=1$.  $H = \langle x, e_0, e_2\rangle$ is abelian and $\alpha = 3 = c-1$.\\
$Y(L) = k[e_2,f]$ where $f = e_2(h^2 + 4xy) + 2(e_0e_1h + e_1^2x - e_0^2y)$,\\
$M = k[x, e_0, e_2, f]$.
\end{itemize}

{\bf \large 3. Preliminaries}\\
\\
{\bf Lemma 3.1.} [BC, p.3]
\begin{itemize}
\item[(1)] Let $H$ be a subalgebra of $L$.  Then $\alpha(H) \leq \alpha(L)$.
\item[(ii)] Suppose $L = L_1 \oplus L_2$ is a direct product.  Let $A_1$ (resp. $A_2$) be an abelian subalgebra of $L_1$ (resp. $L_2$) of maximum dimension.  
Then $A = A_1 \oplus A_2$ is an abelian subalgebra of $L$ of maximum dimension.  In particular,
$$\alpha(L) = \alpha(L_1) + \alpha (L_2)$$
\end{itemize}

{\bf Lemma 3.2.}\\
Any $n$-dimensional Lie algebra $L$ with $i(L) \geq n-2$ satisfies the $R$-property.\\
\ \\
{\bf Proof.} We may assume that $L$ is not abelian, i.e. $i(L) < n$.  Hence $i(L) = n-2$ and $c(L) = (n+n-2)/2 = n-1$. By formula (\ref{3})
$$n-2 = i(L) \leq \alpha(L) \leq c(L) = n-1$$
Consequently, $\alpha(L) \geq n-2 = c(L) - 1$.\hfill $\square$\\
\ \\
{\bf Lemma 3.3.} Let $H$ be a subalgebra of $L$ of codimension one.  Then,
\begin{itemize}
\item[(i)] If $L$ admits a CP $P$ then $H$ satisfies the $R$-property.
\item[(ii)] If $H$ admits a CP $Q$ then $L$ satisfies the $R$-property.\\
\end{itemize}

{\bf Proof.} First we recall from [EO, Proposition 1.6] that we have either $i(H) = i(L) + 1$ (i.e. $c(H) = c(L)$) or $i(H) = i(L) - 1$ (i.e. $c(H) = c(L)-1$).
\begin{itemize} 
\item[(i)] $P$ is a commutative subalgebra of $L$ with $\dim P = c(L)$.  There are two cases to consider:
\begin{itemize}
\item[(1)] $P \subset H$.  Then $c(L) = \dim P \leq c (H) \leq c(L)$.  Therefore $\dim P = c(H)$, i.e. $P$ is a CP of $H$.
\item[(2)] $P\not\subset H$.  Then $P \cap H$ is a commutative subalgebra of $H$ and 
$$c(L) -1 = \dim P-1 = \dim (P \cap H) \leq c(H) \leq c(L)$$
\end{itemize}
There are 2 possibilities:
\begin{itemize}
\item[(2a)] $\dim (P\cap H) = c(H)$, i.e. $P \cap H$ is a CP of $H$.
\item[(2b)] $\dim (P\cap H) \neq c(H)$.  Then $c(H) = c(L)$ and $\dim (P \cap H) = c(L) - 1 = c(H) - 1$.  Hence $H$ satisfies the $R$-property.
\end{itemize}
\item[(ii)]  $Q$ is a commutative subalgebra of $H$ with $\dim Q = c(H) \leq c(L)$.   We have to consider 2 cases:
\begin{itemize}
\item[(1)] $c(H) = c(L)$.  Then $\dim Q = c(L)$ and $Q$ is a CP of $L$
\item[(2)] $c(H) = c(L) - 1$.  Then $\dim Q = c(H) = c(L) - 1$.
\end{itemize}
We may conclude that $L$ satisfies the $R$-property. \hfill $\square$
\end{itemize}

{\bf Remark 3.4.} Let $H$ be a abelian subalgebra of $L$ with $\dim H = c(L) - 1$.
\begin{itemize} 
\item[(i)] If $f \in Y(L)$ but $f \notin S(H)$ then $S(H)k[f]$ is a complete Poisson commutative subalgebra of $S(L)$.
\item[(ii)] If $F(L) \subset H$ then $M = S(H) Y(L)$ is not complete\\
(Indeed, $Y(L) \subset S(F(L)) \subset S(H)$ and so $M = S(H)$).
\end{itemize}

{\bf Proposition 3.5.}\\
Let $L$ be an algebraic Lie algebra without proper semi-invariants.  Assume that $H$ is an abelian subalgebra of dimension $c(L)-1$, which does not contain $F(L)$.  Then $M = S(H)Y(L)$ is a complete,
Poisson commutative subalgebra of $S(L)$.\\
\ \\
{\bf Proof.}  Suppose $Y(L) \subset S(H)$.  Then the second characterization of $F(L)$ (see Section 2) asserts that $F(L) \subset H$, which contradicts our assumption.  Therefore $Y(L) \not\subset S(H)$, 
which implies that $M$ is complete. \hfill $\square$\\
\ \\
The following shows that the condition on the proper semi-invariants cannot be omitted.\\
\ \\
{\bf Example 3.6.}\\
Let $L$ be the algebraic Lie algebra with basis $x, y, z$ and nonzero brackets $[x,y] = y$ and $[x,z] = z$.  Then $i(L) = 1$, $c(L) = 2$, $Y(L) = k$, $F(L) = \langle y, z\rangle$.\\
Now, take $H = \langle x \rangle$; which is an abelian subalgebra of $L$ with $\dim H = 1 = c(L) - 1$.  Clearly $F(L) \not\subset H$.  On the other hand, $M = S(H) Y(L) = k[x]$ is not complete.\\
\ \\
In the following situation the condition that $F(L) \not\subset H$ is automatically satisfied.\\
\ \\
{\bf Proposition 3.7}\\
Assume that $L$ is nonsingular (i.e. $p_L = 1$) without proper semi-invariants.  Then any abelian subalgebra $H$ of $L$ of dimension $c(L) - 1$ does not contain $F(L)$.\\
\ \\
{\bf Proof.}  Suppose $F(L) \subset H$.  Then $\dim F(L)\leq \dim H = c(L) - 1$ and $F(L)$ is abelian. Since $L$ is nonsingular the latter implies that $F(L)$ is a $CP$ of $L$ [O5, Theorem 22].
In particular, $\dim F(L) = c(L)$.  Contradiction.\hfill $\square$\\
\ \\
{\bf Lemma 3.8.}\\
Assume that each $(n-1)$-dimensional, $n \geq 2$, solvable (respectively nilpotent) Lie algebra satisfies the $R$-property.  Then the same holds for any $n$-dimensional solvable (resp. nilpotent) Lie algebra $L$ which is not quasi quadratic.\\
\ \\
{\bf Proof.}  First we see that $L \neq F(L)$ since $F(L)$ is quasi quadratic (as $F(F(L)) = F(L)$).  Being solvable, $L$ admits an ideal $H$ of codimension one containing $F(L)$.
By [EO, Proposition 1.6(4)] $i(H) = i(L) + 1$ which implies that $c(H) = c(L)$.  By assumption $H$ satisfies the $R$-property, i.e. $\alpha(H) \geq c(H) - 1$.  Finally we observe that
$$c(L) - 1 = c(H) - 1  \leq \alpha(H) \leq \alpha (L)$$
\hfill $\square$\\
\ \\
We will now study the occurrence of the $R$-property in various Lie algebras.\\
\ \\
{\bf {\large 4. The nonsolvable case}}\\
First we note that $sl(2,k)$ is the only semi-simple Lie algebra with the $R$-property.\\
\ \\
{\bf Theorem 4.1.}\\
The $R$-property holds for each indecomposable, nonsolvable Lie algebra $L$ of dimension at most seven.  In dimension eight there are 3 exceptions, one of which is $sl(3,k)$.\\
\ \\
{\bf Proof.}  We verify this for each member of the list [O6, pp. 125-136], see also [OAV, pp.554-580].  First we look for a $CP$ in $L$.  If this does not exist it suffices to find an 
abelian subalgebra $H$ with $\dim H = c(L) - 1$.\\
\ \\
{\bf (i) $\boldmath{\dim L \leq 7}$.}\\
$sl(2,k)\ (H = \langle h\rangle, i = 1 =  \alpha = c -1)$ (quadratic),\\
$L_5\ (H = \langle e_0, e_1\rangle, i=1, \alpha = 2 = c-1)$ (quasi quadratic),\\
$L_{6,1}\ (H = \langle e_0, e_1, e_2\rangle, i=2, \alpha = 3 = c-1)$ (quadratic)\\
$L_{6,3}\ (H = \langle x, e_0, e_2\rangle, i=2, \alpha = 3 = c-1)$ (quasi quadratic).  This is Example (5) of Section 2.\\
$L_{6,4}\ (H = \langle e_0, e_1\rangle, i=0, \alpha = 2 = c-1)$ (Frobenius)\\
$L_{7,1}\ (H = \langle e_0, e_1, e_2, e_3\rangle, i=1, \alpha = 4 = c)$\\
$L_{7,2}\ (H = \langle e_0, e_1, e_2, e_3\rangle, i=1, \alpha = 4 = c)$\\
$L_{7,7}\ (H = \langle e_0, e_1, e_2\rangle, i=1, \alpha = 3 = c-1)$\\
$L_{7,8}\ (H = \langle e_0, e_1, e_2\rangle, i=1, \alpha = 3 = c-1)$\\
$L_{7,9}\ (H = \langle x, e_0, e_2\rangle, i=1, \alpha = 3 = c-1)$\\
\ \\
{\bf (ii) $\boldmath{\dim L = 8}$.}\\
$L_{8,1}\ (H = \langle e_0, e_1, e_2, e_3, e_4\rangle, i=2, \alpha = 5 = c)$\\
$L_{8,2}\ (H = \langle e_0, e_1, e_2, e_3, e_4\rangle, i=2, \alpha = 5 = c)$\\
$L_{8,13}\ (H = \langle x, e_0, e_2, e_3\rangle, i=2, \alpha = 4 = c-1)$ (quasi quadratic)\\
$L_{8,14}\ (H = \langle x, e_0, e_2, e_3\rangle, i=2, \alpha = 4 = c-1)$\\
$L_{8,15}\ (H = \langle x, e_0, e_2, e_3\rangle, i=2, \alpha = 4 = c-1)$ (quasi quadratic)\\
$L_{8,17}\ (H = \langle e_0, e_1, e_2, e_3\rangle, i=2, \alpha = 4 = c-1)$ (quasi quadratic)\\
$L_{8,19}\ (H = \langle e_0, e_1, e_2, e_3\rangle, i=0, \alpha = 4 = c)$ (Frobenius)\\
$L_{8,20}\ (a \neq -1)\ (H = \langle e_0, e_1, e_2, e_3\rangle, i=0, \alpha = 4 = c)$ (Frobenius)\\
$L_{8,20}\ (a = -1)\ (H = \langle e_0, e_1, e_2, e_3\rangle, i=2, \alpha = 4 = c-1)$ (quasi quadratic)\\
$L_{8,21}\ (H = \langle e_0, e_1, e_2, e_3\rangle, i=2, \alpha = 4 = c-1)$\\
$L_{8,22}\ (H = \langle e_0, e_1, e_2, e_3\rangle, i=2, \alpha = 4 = c-1)$\\
$L_{8,23}\ (H = \langle x, e_0, e_2, e_3\rangle, i=2, \alpha = 4 = c-1)$\\
$L_{8,25}\ (H = \langle e_0, e_1, e_2, e_3\rangle, i=2, \alpha = 4 = c-1)$\\
$L_{8,26}\ (H = \langle e_0, e_1, e_2, e_3\rangle, i=2, \alpha = 4 = c-1)$\\
$L_{8,27}\ (H = \langle x, e_0, e_2, e_3\rangle, i=2, \alpha = 4 = c-1)$\\
$L_{8,28}\ (H = \langle e_0, e_1, e_2, e_3\rangle, i=0, \alpha = 4 = c)$ (Frobenius)\\
$L_{8,24} = sl(3,k) = A_2$, $i=2$, $c=5$.  From the first row of Table 1 we see that $\alpha = \left[\frac{9}{4}\right] = 2 < 4=c-1$.  So, $sl(3,k)$ does not have the $R$-property.\\
The same is true for the following.  Indeed,\\
$L_{8,16}$ ($H = \langle e_0, e_1, e_2\rangle$, $i=2$, $c=5$, $\alpha = 3 < c - 1$) (quasi-quadratic)\\
$L_{8,18}$ ($H = \langle e_0, e_1, e_2\rangle$, $i=2$, $c=5$, $\alpha = 3 < c - 1$) (quasi-quadratic) \hfill $\square$\\
\ \\
In [EO] we studied the occurrence of CP's in various Lie algebras.  In particular we have\\
{\bf Proposition 4.2} [EO, p.142]\\
Let $\mathfrak{g}$ be simple and $V$ an irreducible $\mathfrak{g}$-module with $\dim \mathfrak{g} < \dim V$.  Consider the semi-direct product $L = \mathfrak{g} \oplus V$.  
Then $V$ is a CP of $L$.  So $L$ satisfies the $R$-property.\\
\ \\
{\bf {\large 5. The solvable case}}\\
{\bf Lemma 5.1.} Let $L$ be an $n$-dimensional solvable Lie algebra with minimal index (i.e. $i(L) = 0$ (or $1$) if $n$ is even (or odd)).  If $n 
\leq 10$, $n \neq 9$ then $L$ satisfies the $R$-property.\\
Moreover, $L$ admits a CP if in addition $L$ is Frobenius with $n \leq 6$.\\
\ \\
{\bf Proof.} Since $L$ is solvable we obtain from the Introduction that
$$]\frac{1}{2} (\sqrt{8n+9} - 3)[\ \leq \alpha (L) \leq c(L)$$
where we denote by $]x[$ the least integer greater than or equal to the real number $x$.\\
Now, from Table 2 below we deduce that 
$$\alpha(L) = c(L)\ \mbox{if}\ i(L) = 0\ \mbox{(i.e. $L$ is Frobenius) and}\ n \leq 6$$
(which means that $L$ contains a CP) and also that:
$$\alpha(L) \geq c(L) - 1\ \mbox{if}\ n \leq 10, n \neq 9$$
in other words the $R$-property is satisfied. \hfill $\square$

\begin{center}
{\bf Table 2}\\
\ \\
\begin{tabular}{l|ccccccccc}
$n$ &2 &3 &4 &5 &6 &7 &8 &9 &10\\
\hline
$i(L)$ &0 &1 &0 &1 &0 &1 &0 &1 &0\\
$c(L)$ &1 &2 &2 &3 &3 &4 &4 &5 &5\\
$]\frac{1}{2} (\sqrt{8n+9}-3)[$ &1 &2 &2 &2 &3 &3 &3 &3 &4
\end{tabular}
\end{center}

{\bf Theorem 5.2.} Let $L$ be solvable of dimension $n$.
\begin{itemize}
\item[(1)] $L$ has the $R$-property if $n \leq 6$.
\item[(2)] Let $n = 7$.  Then $L$ has the $R$-property if one of the following conditions is satisfied:
\begin{itemize}
\item[(i)] $i(L) \neq 3$
\item[(ii)] $i(L) = 3$ and $L$ is not quasi quadratic
\end{itemize}
\end{itemize}

{\bf Proof.} We may assume that $L$ is not abelian.
\begin{itemize}
\item[(1)] First we suppose that $n \leq 5$.  Then the result follows from Lemma 5.1 and from Lemma 3.2 since $i(L)$ is either minimal or is equal to $n-2$.\\
Next, take $n = 6$.\\
Then the case $i(L) = 0$ or $i(L) = 4$ can be treated in the same way as above.  In the remaining case that $i(L) = 2$ we observe that
$$3 = ] \frac{1}{2} (\sqrt{8n+9} - 3) [\ \leq \alpha (L) \leq c(L) = 4$$
i.e. $\alpha(L) \geq 3 = c(L) - 1$.
\item[(2)] $n=7$
\begin{itemize}
\item[$\bullet$] If $i(L) \neq 3$ then the result follows at once from Lemma 5.1 and from Lemma 3.2 since we have either $i(L) = 1$ or $i(L) = 5 = n-2$.
\item[$\bullet$] If $L$ is not quasi quadratic, then it suffices to combine (1) with Lemma 3.8 \hfill $\square$
\end{itemize}
\end{itemize}
{\bf Corollary 5.3.}\\
Any 7-dimensional solvable Lie algebra with trivial center satisfies the $R$-property.\\
\ \\
The following shows that condition (ii) of Theorem 5.2 cannot be removed.\\
\ \\
{\bf Example 5.4.}\\
Let $L$ be the solvable Lie algebra with basis $x_1,\ldots, x_7$ and nonzero brackets:\\
$[x_1,x_3] = -x_3$, $[x_1, x_6] = x_6$, $[x_2,x_4] = -x_4$, $[x_2,x_5] = x_5$, $[x_3, x_6] = x_7$, $[x_4,x_5] = x_7$.\\
$L$ is quasi quadratic (but not quadratic) and $H = \langle x_5, x_6, x_7\rangle$ is abelian.  Also $\alpha(L) = 3 = i(L)$ and $c(L) = 5$, $p_L = 1$.  So, $L$ does not have the $R$-property and Proposition 3.5 is not applicable.
However, by [JS, 5.7], [O6, Theorem 29]
$$Y(L) = k[x_7, f, g]\ \mbox{where}\ f= x_1x_7 - x_3x_6\ \mbox{and}\ g = x_2x_7 - x_4x_5$$
and $M = S(H) Y(L) = k[x_5, x_6, x_7, f,g]$ is still a complete Poisson commutative subalgebra of $S(L)$.  Moreover, it clearly satisfies Milovanov's conjecture.  We now want to demonstrate that $M$ is also
strongly complete.  Take $\xi \in L^\ast$.  Using the fact that $d_{\xi}(x) = x$ for all $x \in L$ we see that
\begin{eqnarray*}
\begin{array}{l}
d_\xi(f) = \xi(x_1)x_7 + x_1\xi(x_7) - \xi(x_3)x_6 - x_3\xi (x_6)\\
d_\xi(g) = \xi(x_2)x_7 + x_2\xi(x_7) - \xi(x_4)x_5 - x_4\xi (x_5)
\end{array}
\end{eqnarray*}
Next we consider the Jacobian locus of the generators of $M$:
\begin{eqnarray*}
J &=& \{\xi \in L^\ast \mid x_5, x_6, x_7, d_\xi (f), d_\xi(g)\ \mbox{are linearly dependent}\}\\
&=& \{ \xi \in L^\ast \mid \xi (x_5) = \xi(x_6) = \xi(x_7) = 0\}
\end{eqnarray*}
Clearly, $\codim \ J = 3 \geq 2$.  By combining [PPY, Theorem 1.1] and [PY, 2.1] we may conclude that $M$ is strongly complete. \hfill $\square$\\
\ \\
{\bf Example 5.5.}\\
Let $B$ be the Borel subalgebra of a simple Lie algebra of type $B_3$.  We know that $\dim B = 12$ and $i(B) = 0$ (so $B$ is Frobenius), 
$c(B) = 6$ and $\alpha(B) = 5$ by Table 1, see also [EO, p.146; Su].  Hence, $B$ does not have a CP, but the $R$-property holds.\\
\ \\
{\bf Example 5.6. The generalized diamond Lie algebra $\boldmath{D_n}$ ($\boldmath{n}$ even)}\\
Let $H$ be the $(2m+1)$-dimensional Heisenberg Lie algebra with canonical basis $x_1,\ldots, x_m, y_1,\ldots, y_m, z$ and nonzero brackets: $[x_i, y_i] = z$, $i=1,\ldots, m$.\\
Next consider the semi-direct product $D_n = kt \oplus H$, where $t$ is the derivation of $H$, defined by
$$t(x_i) = x_i, t(y_i) = -y_i, i=1,\ldots, m\ \mbox{and}\ t(z) = 0$$
and where $n = 2(m+1)$.\\
Then $D_n$ is solvable and quadratic.  Indeed $D_n$ admits a nondegenerate, invariant, symmetric bilinear form $b$ given by the following nonzero 
entries $b(t,z) = 1$, $b(x_i,y_i) = 1$, $i = 1,\ldots, m$\\
One verifies that $i(D_n) = 2$ and $c(D_n) = m+2$.  By [O6, Theorem 52]
$$Y(D_n) = k[z,f],\ \mbox{where}\ f = tz + x_1y_1 + \ldots + x_m y_m$$
$f$ is the Casimir of $D_n$ w.r.t. $b$.
$D_n$ being quadratic, is also quasi quadratic and so does not possess any CP's, i.e. $\alpha(D_n) \leq c(D_n) - 1$.  On the other hand, $A = \langle y_1,\ldots, y_m,z\rangle$ is abelian with $\dim A = m+1 = c(D_n) - 1$.
Consequently $\alpha (D_n) = c(D_n) - 1 = m+1 = \frac{n}{2}$ and $D_n$ has the $R$-property.\\
Finally, $M = S(A) Y(D_n) = k [y_1,\ldots, y_m, z,f]$ is a strongly complete Poisson commutative subalgebra of $S(D_n)$, which satisfies the Milovanov conjecture.\\
\ \\
{\bf {\large 6. The nilpotent case}}\\
{\bf Lemma 6.1.} Let $L$ be a $n$-dimensional nilpotent Lie algebra with minimal index (i.e. $i(L) = 1$ if $n$ is odd and $i(L) = 2$ if $n$ is even).  If $n \leq 11$, $n \neq 10$ then $L$ satisfies the $R$-property.
Moreover, $L$ admits a CP if $n \leq 7$, $n \neq 6$.\\
\ \\
{\bf Proof.} Since $L$ is nilpotent we know from the Introduction that
$$]\frac{1}{2}(\sqrt{8n+1} - 1) [\ \leq \alpha (L) \leq c(L)$$
Now, from Table 3 below we deduce that:
$$\alpha (L) = c(L)\ \mbox{if}\ n\leq 7\ \mbox{and}\ n \neq 6$$
which means that $L$ admits a CP.\
We also notice that
$$\alpha(L) \geq c(L) -1\ \mbox{if}\ n \leq 11, n \neq 10$$
which implies that $L$ satisfies the $R$-property.\hfill $\square$
\begin{center}
{\bf Table 3}\\
\ \\
\begin{tabular}{l|cccccccccc}
$n$ &2 &3 &4 &5 &6 &7 &8 &9 &10 &11\\
\hline
$i(L)$ &2 &1 &2 &1 &2 &1 &2 &1 &2 &1\\
$c(L)$ &2 &2 &3 &3 &4 &4 &5 &5 &6 &6\\
$]\frac{1}{2} (\sqrt{8n+1}-1)[$ &2 &2 &3 &3 &3 &4 &4 &4 &4 &5
\end{tabular}
\end{center}

{\bf Remark.} There exist 6-dimensional nilpotent Lie algebras of index 2 without CP's.  For instance $\mathfrak{g}_{6,18}$ (\# 21 of [O4]) and $\mathfrak{g}_{6,20}$ (\# 28 of [O4]).\\
\ \\
In (1) of the following we give a proof for Rentschler's case by case observation mentioned in the Introduction.\\
\ \\
{\bf Theorem 6.2.} Let $L$ be nilpotent of dimension $n$.
\begin{itemize}
\item[(1)] If $n\ \leq 7$ then $L$ satisfies the R-property\\
More precisely:
\begin{itemize}
\item[(a)] $\alpha(L) = c(L)$ if $F(L)$ is abelian.\\
If in addition $L$ is nonsingular then $F(L)$ is the only CP of $L$.
\item[(b)] $\alpha(L) = c(L) - 1$ if $F(L)$ is not abelian.
\end{itemize}
\item[(2)] Now assume that $n=8$.  Then $L$ satisfies the R-property if one of the following holds:
\begin{itemize}
\item[(i)] $i(L) \neq 4$.
\item[(ii)] $i(L) = 4$ and $L$ is not quasi quadratic.
\end{itemize}
\end{itemize}

{\bf Proof.} 
\begin{itemize}
\item[(1)] By Theorem 5.2 and Lemma 3.8 it suffices to show that any 7-dimensional quasi quadratic nilpotent Lie algebra $L$ satisfies the R-property.  We may assume that $L$ is not abelian.  Then $L$ does not possess
a CP (i.e. $\alpha(L) \leq c(L)-1$) and $i(L) \geq \dim Z(L) \geq 2$ (see Section 2).  Consequently, we have either $i(L) = 5$ or $i(L) = 3$.
\begin{itemize}
\item[(i)] If $i(L) = 5$ we are done by Lemma 3.2 and in this situation $\alpha(L) = 5 = i(L)$.
\item[(ii)] Now suppose $i(L) = 3$.  Then $c(L) = 5$ and so we need to construct a 4-dimensional abelian subalgebra $H$ of $L$.  This is easy if $\dim Z(L) = 3$ (simply take $H= Z(L) \oplus kw$ with 
$w\in L\backslash Z(L)$).\\
So we may assume that $\dim Z(L) = 2$, say $Z(L) = \langle y,z\rangle$.  Let $Z_1(L)$ be the next ideal of the upper central series i.e.
$$Z_1(L) = \{x \in L \mid [L,x] \subset Z(L)\}$$
Then $Z(L) \subs\limits_{\neq} Z_1(L)$ as $L$ is nilpotent.  Choose $v \in Z_1(L) \backslash Z(L)$ and consider a subspace $U$ of $L$ such that
$$U \oplus \langle v,y,z\rangle = L$$
Then there are $\lambda, \mu \in U^\ast$ such that
$$[u,v] = \lambda (u)y + \mu(u)z\ \mbox{for all}\ u \in U$$
Since $\dim U = 4$ we can find a nonzero $u \in U$ such that $\lambda(u) = 0 = \mu(u)$, i.e. $[u,v] = 0$.\\
Finally $H = \langle u,v,y,z\rangle$ is a 4-dimensional abelian subalgebra of $L$.  Note that $\alpha(L) = 4 = c(L) - 1$.  The remainder of (i) follows directly from [O5, p.93].
\end{itemize}
\item[(2)] 
\begin{itemize}
\item[(i)] $i(L) \neq 4$.  We may assume that $L$ is not abelian.  Hence either $i(L) = 2$ or $i(L) = 6$ since $n = 8$.  In the first case it suffices to apply Lemma 6.1 and in the second case Lemma 3.2.
\item[(ii)]  Now the result follows from (1) combined with Lemma 3.8.\hfill $\square$
\end{itemize}
\end{itemize}

{\bf Corollary 6.3} Let $L$ be an 8-dimensional nilpotent Lie algebra with $\dim Z(L) = 1$.  Then $L$ satisfies the $R$-property.\\
\ \\
The following example shows that condition (ii) in Theorem 6.2 cannot be removed.\\

{\bf Example 6.4.} Let $L$ be an 8-dimensional nilpotent Lie algebra with basis $x_1,\ldots, x_8$ and nonzero brackets:\\
$[x_1,x_3] = x_6$, $[x_1,x_4] = x_5$, $[x_1,x_5] = x_7$, $[x_2, x_3] =x_5$,\\
$[x_2,x_4] = x_6$, $[x_2,x_6] = x_7$, $[x_3,x_5] = x_8$, $[x_4, x_6] =x_8$.\\
$L$ is quasi quadratic (but not quadratic).\\
$H = \langle x_5, x_6, x_7, x_8\rangle$ is abelian.  $\alpha(L) = 4 = i(L)$ and $c(L) = 6$.  Hence the R-property is not valid here.\\
$Y(L) = k[x_7,x_8,f,g]$ where $f = x_1x_8 - x_3x_7 + x_5x_6$, $g = x_5^2 + 2x_2x_8 - 2x_4x_7 + x_6^2$.\\
$M = S(H) Y(L) = k[x_5, x_6, x_7, x_8, x_1x_8 - x_3x_7, x_2x_8 - x_4x_7]$, which is strongly complete.  Clearly $L$ satisfies Milovanov's conjecture.\\
\ \\
{\bf Example 6.5.} Let $L$ be the 8-dimensional Lie algebra with basis $x_1,\ldots,x_8$ and nonzero brackets:\\
$[x_1,x_2] = x_5$, $[x_1,x_3] = x_6$, $[x_1,x_4] = x_7$, $[x_1, x_5] =-x_8$, $[x_2,x_3] = x_8$,\\
$[x_2,x_4] = x_6$, $[x_2,x_6] = -x_7$, $[x_3,x_4] = -x_5$, $[x_3, x_5] = -x_7$, $[x_4,x_6] = -x_8$.\\
$L$ is characteristically nilpotent [DL] of index 2, with center $Z(L) = \langle x_7, x_8\rangle$.  Since $i(L) = 2 = \dim Z(L)$ $L$ is square integrable
and $F(L) = Z(L)$.  This implies that $Y(L) = S(Z(L)) = k [x_7, x_8]$. By [EO, Remark 4.2(b)] $L$ does not have any CP's, i.e. $\alpha(L) < c(L) = 5$.  On the other hand, 
$H = \langle x_5, x_6, x_7, x_8\rangle = [L,L]$ is abelian and
$\dim H = 4 = c(L) - 1$.  So, the R-property holds and $\alpha(L) = 4$.  Note that $F(L) = Z(L)$ is contained in $H$.  Therefore $S(H) Y(L) = k [x_5, x_6, x_7, x_8]$ is not complete.  However put
$f = x_1x_7 + x_2x_8 - x_3x_8 - x_4x_7$.  Then $M = k [x_5, x_6, x_7, x_8, f]$ is a strongly complete, Poisson commutative subalgebra of $S(L)$.  Clearly, $L$ satisfies tha Milovanov conjecture.\\
\ \\
{\bf Theorem 6.6.} [EO, Theorem 6.2]\\
Let $L$ be a simple Lie algebra of type $A$ or $C$, $P$ a parabolic subalgebra of $L$.  Then the nilradical of $P$ admits a CP.\\
\ \\
We can now easily rediscover the value of $\alpha(L)$ for some standard Lie algebras.  As a bonus we will obtain a short proof for Schur's formula (Corollary 6.10).\\
\ \\
{\bf Proposition 6.7.}\\
Let $H_n$ be the standard $n$-dimensional Heisenberg Lie algebra with basis\\
$x_1, \ldots, x_m, y_1,\ldots, y_m, z$ (so $n = 2m+1$) with nonzero brackets $[x_i,y_i] = z$, $i = 1,\ldots, m$.  Then
$$P = \langle y_1,\ldots, y_m,z\rangle$$
is a CP of $H_n$ and therefore $\alpha (H_n) = m+1$.\\
\ \\
{\bf Proof.} It is easy to see that $i(H_n) = 1$ and $c(H_n) = (2m+1+1)/2 = m+1$.  On the other hand $P$ is a commutative ideal of $H_n$ of dimension $m+1$ and hence is a CP of $H_n$. \hfill $\square$\\
\ \\
{\bf Lemma 6.8.}
\begin{itemize}
\item[(1)] Let $N$ be the nilradical of a solvable Lie algebra $L$.  Then $\alpha(N) = \alpha (L$).
\item[(2)] Let $B$ be a Borel subalgebra of an arbitrary Lie algebra $L$.  Then $\alpha (B) = \alpha (L)$.
\end{itemize}

{\bf Proof.} 
\begin{itemize}
\item[(i)] $\alpha(N) \leq \alpha(L)$ since $N\subset L$.  Being solvable, $L$ admits an abelian ideal $A$ such that $\dim A = \alpha(L)$ [BC, Proposition 2.6].  But $N$ is the greatest nilpotent ideal of $L$.  
Therefore $A \subset N$ and so $\alpha(L) = \dim A \leq \alpha (N)$.  Consequently, $\alpha(N) = \alpha(L)$.
\item[(2)] $\alpha(B) \leq \alpha(L)$ since $B \subset L$.  Next, we take an abelian subalgebra $A$ of $L$ such that $\dim A = \alpha(L)$.  In particular, $A$ is a solvable subalgebra of $L$ and therefore it is contained in a Borel subalgebra $B_1$ of $L$.
Hence
$$\alpha(L) = \dim A \leq \alpha (B_1) = \alpha(B)$$
The latter equality is valid because $B$ and $B_1$ are isomorphic [TY, Theorem 29.4.7].  We may conclude that $\alpha(B) = \alpha(L)$.\hfill $\square$
\end{itemize}

{\bf Theorem 6.9.}\\
Let $N_n$ be the Lie algebra of all $n \times n$ strictly upper triangular matrices, which we consider as the nilradical of the standard Borel subalgebra $B_n$ of $sl(n,k)$.  Then,
\begin{itemize}
\item[(i)] $N_n$ admits a CP-ideal $P$ and $\alpha(N_n) = [n^2/4]$, (which is clearly equal to $q^2$ if $n = 2q$ and equal to $q(q+1)$ if $n = 2q + 1)$
\item[(ii)] $\alpha(N_n) = \alpha(B_n) = \alpha (sl(n,k)) = [n^2/4]$
\end{itemize}
In particular $P$ is also an abelian subalgebra of $B_n$ (and of $sl(n,k)$) of maximum dimension.  Note that the latter equality proves the first row of Table 1.\\
\ \\
{\bf Proof.}
\begin{itemize}
\item[(i)] $N_n$ is spanned by the standard matrices $E_{ij}$, $i \leq i < j \leq n$ and $B_n$ is spanned by the same $E_{ij}$'s together with $E_{ii} - E_{jj}$, $1 \leq i < j \leq n$.  
Put $q = [n/2]$ and consider the following abelian ideal of $N_n$ (see [O3, p.285]):
$$P = \langle E_{ij} \mid 1 \leq i \leq q, q +1 \leq j \leq n\rangle$$
Clearly, $\dim P = q(n-q)$.  We claim that $P$ is a CP of $N_n$ (so this is a special case of Theorem 6.6).  We only have to check that $\dim P = c(N_n)$.  It is well known that $i(N_n) = q$ [O3, Theorem 4.1].
\begin{center}
Hence, $c(N_n) = (\dim N_n + i(N_n))/2 = (\frac{1}{2} n(n-1) + q)/2$
\end{center}
So, if $n = 2q$ then
$$c(N_n) = (q(2q-1) + q)/2 = q^2 = q(n-q) = \dim P$$
On the other hand, if $n = 2q+1$
$$c(N_n) = ((2q+1)q + q)/2 = q^2 + q = q(q+1) = q(n-q) = \dim P$$
We may conclude that $P$ is indeed a CP of $N_n$ and also that $\alpha(N_n) = c(N_n) = [n^2/4]$ by the above.
\item[(ii)] This now follows directly from (i) and Lemma 6.8. \hfill $\square$
\end{itemize}

{\bf Corollary 6.10.}\\
Let $T_n$ be the Lie algebra of all $n \times n$ upper triangular matrices with coefficients in $k$.  Then $\alpha(T_n) = [n^2/4] + 1 = \alpha(gl(n,k))$ and $kI_n \oplus P$ is an abelian subalgebra of maximum
dimension of both $T_n$ and $gl(n,k)$.\\
\ \\
{\bf Proof.} Clearly, $T_n = kI_n \oplus B_n$ and $gl(n,k) = kI_n \oplus sl(n,k)$ are direct products.  Hence by (ii) of Lemma 3.1 $kI_n \oplus P$ is an abelian subalgebra of maximum dimension of both of them and
$$\alpha(T_n) = 1 +\alpha(B_n) = 1 + [n^2/4] = 1 + \alpha(sl(n,k)) = \alpha (gl(n,k))$$ \hfill $\square$\\
\ \\
{\bf {\large 7. The metabelian case}}\\
The following is a slightly rephrased version of [Mil3, Theorem 2].\\
\ \\
{\bf Theorem 7.1.}\\
Let $L$ be an $n$-dimensional metabelian Lie algebra (i.e. $[L,L] \subset Z(L)$) with $t = \dim [L,L] \geq 2$.  Then $L$ contains an abelian subalgebra of dimension $s = [(2n + t^2 + t)/(t+2)]$ (i.e. $\alpha(L) \geq s$).  
If $L$ is generic in the sense of [Mil3, Definition 1] then $L$ does not contain any $(s+1)$-dimensional abelian subalgebra (i.e. $\alpha(L) =s$).\\
\ \\
{\bf Corollary 7.2.}\\
Let $L$ be a metabelian Lie algebra with $t = \dim [L,L] = 2$ and $i(L) = 2$ or 3.  Then $L$ admits a CP.\\
\ \\
{\bf Proof.}
\begin{itemize}
\item[(1)] First suppose $i(L) = 2$.  Then $n = \dim L$ is even, i.e. $n = 2q$ for some integer $q$.  Hence, $c(L) = (\dim L + i(L))/2 = (2q +2)/2 = q+1$.  By the previous theorem
$$\alpha(L) \geq [(2n + t^2 + t)/ (t+2)] = [(4q + 4 + 2)/4] = [q+1 + \frac{1}{2}] = q + 1 = c(L)$$
By formula (3) $\alpha(L) = c(L)$ which implies that $L$ has a CP.  Note that in this situation $L$ is square integrable (i.e. $\dim Z(L) = i(L)$) and $[L,L] = Z(L)$.  Indeed,
$2 = \dim [L,L] \leq \dim Z(L) \leq i(L) = 2$.
\item[(2)] Next assume that $i(L) = 3$.  Then $\dim L$ is odd, i.e. $\dim L = 2q+1$ for some integer $q$.
Hence $c(L) = (\dim L + i(L))/2 = (2q + 1 + 3)/2 = q+2$.\\
By the previous theorem $\alpha(L) \geq [(2n + t^2 + t)/ (t+2)] = [(4q + 2 + 4 + 2)/4] = [q+2] = q+2 = c(L)$.  Consequently, $\alpha(L) = c(L)$ and so $L$ contains a CP.\hfill $\square$
\end{itemize}

{\bf Example 7.3.} [YD] Consider the 8-dimensional metabolism Lie algebra $L$, with basis $x_1,\ldots, x_8$ and nonzero brackets $[x_1, x_2] = x_7$, $[x_2,x_3] = x_8$, $[x_3,x_4] = x_7$, $[x_4,x_5] = x_8$, $[x_5, x_6] = x_7$.\\
Clearly, $L$ is metabelian since $[L,L] = Z(L) = \langle x_7, x_8\rangle$, which is 2-dimensional. $i(L) = 2$ and $c(L) = 5$.\\
So, $L$ contains a CP, namely $\langle x_1, x_3, x_5, x_7, x_8\rangle$.\\
\ \\
{\bf Theorem 7.4.}\\
Let $L$ be an $n$-dimensional metabelian Lie algebra with $t = \dim [L,L] \geq 2$.  Then $L$ satisfies the R-property if $n \leq 9$.\\
\ \\
{\bf Proof.} We know that
$$2 \leq t \leq i(\mathfrak{g} )(\ast)\ \mbox{and}\ s \leq \alpha(L) \leq c(L)$$
where $s = [(2n+t^2+t)/(t+2)]$ by Theorem 7.1.\\
Clearly it suffices to show that $c(L) - s \leq 1$.  By (1) of Theorem 6.2 we only need to consider the cases $n = 8$ and $n=9$.
\begin{itemize}
\item[(1)] $n = 8$.  Then $i(L)$ is even and we may assume that $2 < i(L)$ (by Corollary 7.2) and also that $i(L) < n-2=6$ (by Lemma 3.2).  Therefore $i(L) = 4$ and $c(L) = (8+4)/2 = 6$.\\
By $(\ast)$ there remain inly 3 cases to examine:\\
(1a) $t = 2$. Then $c(L) - s = 6 - [22/4] = 6-5 = 1$\\
(1b) $t =3$.  Then $c(L) - s = 6 - [28/5] = 6-5 =1$\\
(1c) $t = 4$.  Then $c(L) - s = 6 - [36/6] = 6-6 = 0$
\item[(2)] $n = 9$.  Then $i(L)$ is odd and $i(L) \geq 3$ (by $(\ast)$) and also that $i(L) < n - 2 = 7$ (by Lemma 3.2).\\
So $i(L)$ is 3 or 5 and we have to treat the following cases:\\
(2a) $t = 2$, $i(L) = 3$.  Now use corollary 7.2.\\
(2b) $t= 2$, $i(L) = 5$.  Then
$$c(L) - s = 7 = [24/4] = 7-6=1$$
(2c) $t= 3$, $i(L) = 3$.  Then
$$c(L) - s = 6 - [30/5] = 6-6=0$$
(2d) $t = 3$, $i(L) = 5$. Then
$$c(L) - s = 7 - [30/5] = 7-6=1$$
(2e) $t=5$, $i(L) = 5$.  Then
$$c(L) - s =7 - [48/7] = 7-6 = 1.$$
\hfill $\square$
\end{itemize}

{\bf {\large 8. The filiform case}}\\
{\bf Definition 8.1} Consider the descending central series of $L$
$$C^1(L) = L, C^2(L) = [L,L],\ldots, C^i(L) = [L,C^{i-1}(L)],\ldots$$
which satisfies $[C^i(L), C^j(L)] \subset C^{i+j}(L)$, $i,j \geq 1$.\\
An $n$-dimensional Lie algebra $L$ is called filiform if $\dim C^i(L) = n-i$, $i=2,\ldots,n$.\\
In particular $C^n(L) = 0$ (and thus $L$ is nilpotent) and $Z(L) = C^{n-1}(L)$ is 1-dimensional.\\
\newpage
{\bf Proposition 8.2} [V, p.92], [B, p.24]\\
For any filiform Lie algebra $L$ there exists a so called adapted basis $x_1, x_2,\ldots, x_n$ with the following brackets, the undefined brackets being zero:\\
$[x_1,x_i] = x_{i+1}, i = 2, \ldots, n-1$\\
$[x_i, x_j] \in \langle x_{i+j},\ldots, x_n\rangle\ i,j \geq 2, i + j \leq n$\\
$[x_{i+1}, x_{n-i}] = (-i)^i ax_n, 1 \leq i \leq n-1$\\
with a certain $a \in k$, which is zero if $n$ is odd.\\
Moreover, the brackets $[x_i,x_j]$ for $i,j \geq 2$ are completely determined by the brackets
$$[x_i,x_{i+1}] = \sum\limits_{j=2i+1}^n a_{ij} x_j,\ 2 \leq i \leq [n/2]$$
It is also easy to see that\\
$C^2(L) = \langle x_3, x_4, \ldots, x_n\rangle,\ldots, C^i (L) = \langle x_{i+1},\ldots, x_n\rangle ,\ldots, C^{n-1} (L) = \langle x_n\rangle$.\\
\ \\
{\bf Example 8.3.} The standard filiform Lie algebra $L_n$.\\
This is the filiform Lie algebra with basis $x_1,\ldots, x_n$, $n\geq 3$ with nonzero brackets $[x_1, x_i] = x_{i+1}$, $i = 2,\ldots, n-1$.\\
$i(L_n) = n-2$ and $c(L_n) = n-1$.  Clearly $P = \langle x_2, x_3,\ldots, x_n\rangle$ is a CP, which coincides with $F(L_n)$ if $n \geq 4$.  Moreover $L_n$ is coregular if and only if $n \leq 4$ [OV, Example 1.7], [O6, Theorem 51].\\
\ \\
{\bf Proposition 8.4.} See [C, Proposition 3.2], [BC, Proposition 5.6]\\
Let $L$ be an $n$-dimensional nonstandard filiform Lie algebra.  Let $m$ be the smallest integer such that $C^m(L)$ is abelian.  Then $C^m(L)$ is the unique abelian ideal of $L$ of maximum dimension.\\
In particular, $\alpha(L) = \dim C^m(L) = n-m$.  Furthermore,
$$m = \max \{i \mid [x_i,x_{i+1}]\neq 0\}$$
where $x_1,\ldots, x_n$ is an adapted basis of $L$.\\
\ \\
{\bf Proof.} (of the last statement only)\\
Consider $C^{m-1}(L)$ with its basis $x_m, x_{m+1},\ldots, x_n$.  As $C^m(L)$ is abelian, the structure matrix $M = ([x_i,x_j])$, $i,j=m,\ldots, n$ of $C^{m-1}(L)$ is given by:
\begin{eqnarray*}
\begin{array}{c|cccc}
       &x_m &x_{m+1} &\ldots &x_n\\
 \hline
 x_m &0 &[x_m, x_{m+1}] &\ldots &[x_m,x_n]\\
 x_{m+1} &-[x_m,x_{m+1}] &0 &\ldots &0\\
       &\vdots &\vdots & &\vdots\\
 x_n &-[x_m, x_n] &0 &\ldots &0
\end{array}
\end{eqnarray*}
Since $[x_s,x_{s+1}]=0$, $s \geq m+1$, it suffices to show that $[x_m, x_{m+1}] \neq 0$.  So, let us suppose that $[x_m, x_{m+1}]=0$.  Then we can show that $[x_m, x_{m+r}] = 0$, $r = 1,\ldots, n-m$ by induction on $r$.
This is clear if $r=1$.  Next take $r \geq 2$ and assume that $[x_m, x_{m+r}] = 0$.\\
Then
\begin{eqnarray*}
0 &=& [x_1, [x_m, x_{m+r}]] = [[x_1,x_m], x_{m+r}] + [x_m, [x_1, x_{m+r}]]\\
  &=& [x_{m+1}, x_{m+r}] + [x_m, x_{m+r+1}] = [x_m, x_{m+r+1}]
\end{eqnarray*}
the first term being zero because $C^m(L)$ is abelian.  Consequently $M = 0$, i.e. $C^{m-1}(L)$ is abelian, which contradicts the assumption of the proposition. \hfill $\square$\\
\ \\
{\bf Proposition 8.5.}\\
Let $L$ be an $n$-dimensional filiform Lie algebra.  Put $q = [(n+1)/2]$.  Then $H = C^{n-q}(L)$ is a $q$-dimensional abelian ideal of $L$, i.e. $\alpha(L) \geq q$.\\
Assume in addition that $L$ has minimal index , then $L$ satisfies the R-property.  More precisely:
\begin{itemize}
\item[(1)] If $i(L) = 1$ (so $n$ is odd) then $H$ is a CP of $L$ and $\alpha(L) = q$.
\item[(2)] If $i(L) = 2$ (so $n$ is even) then $\dim H = q = c(L)-1$ (but the existence of a CP is still possible, see (3b) of Theorem 8.6).
\end{itemize}

{\bf Proof.} Let $x_1,\ldots, x_{q-1}, x_q, x_{q+1},\ldots, x_n$ be an adapted basis of $L$. We have to consider 2 cases
\begin{itemize}
\item[(1)] $n$ is odd. Then $n = 2q-1$ and $n-q = q-1$.  It turns out that $H = C^{q-1}(L) = \langle x_q, x_{q+1},\ldots, x_n\rangle$ is a $q$-dimensional abelian ideal.\\
Next we assume in addition that $i(L) = 1$.  Then
$$c(L) = (n+1)/2 = (2q-1+1)/2 = q= \dim H$$
Therefore $H$ is a CP of $L$.
\item[(2)] $n$ is even.  Then $n = 2q$ and
$$H = C^{n-q}(L) = C^q(L) = \langle x_{q+1},\ldots,x_n\rangle$$
which is a $q$-dimensional abelian ideal of $L$ since
$$[H,H] = [C^q(L), C^q(L)] \subset C^{2q}(L) = C^n(L) = 0$$
Now assume in addition that $i(L) = 2$.  Then $c(L) = (2q+2)/2 = q+1$.  Hence, $\dim H = q = c(L) - 1$.\hfill $\square$
\end{itemize}

{\bf Theorem 8.6.}\\
Let $L$ be filiform of one of the major types $L_n$, $Q_n$, $R_n$ and $W_n$ [GK, p.111], [O6, pp.120-121].  Then the R-property holds.  In fact $L_n$, $R_n$, $W_n$ admit a CP.
Furthermore, $\alpha(L_n) = n-1$, $\alpha(Q_n) = n/2$, $\alpha(R_n) = n-2$, $\alpha(W_n) = [(n+2)/2]$.  Also, the Milovanov conjecture is valid for $L$.\\
\ \\
{\bf Proof.}
\begin{itemize}
\item[(1)] The case where $L$ is type $L_n$ has already been done in Example 8.3.
\item[(2)] Suppose $L$ is of type $Q_n$.\\
Basis: $x_1,\ldots,x_n$, $n = 2q$\\
Nonzero brackets: $[x_1, x_i] = x_{i+1}$, $i = 2,\ldots, n-2$\\
and $[x_j, x_{n-j+1}] = (-1)^{j+1} x_n$, $j = 2,\ldots, q$.\\
Note that our basis differs slightly from an adapted basis.  Put $\xi = x_n^\ast \in L^\ast$.  it is easy to see that $\xi$ is regular and that $L(\xi) = \langle x_1, x_n\rangle$.\\
Then $i(L) = \dim L (\xi) = 2$ and $c(L) = (n+2)/2 = q+1$.\\
By Proposition 8.5 $L$ satisfies the $R$-property.  More precisely, $H = C^q(L) = \langle x_{q+1},\ldots x_n\rangle$ is an abelian ideal of $L$ and $\dim H = q =c(L) - 1$. \\
So $\alpha(L) \geq c(L)-1$.  Furthermore,
$$F(L) = \langle x_1, x_3,\ldots, x_n\rangle$$
Indeed, $[L,L(\xi)] \subset [L, F(L)] \subset F(L)$.  Hence $x_1, x_3, \ldots, x_n \in F(L)$.  On the other hand, the centralizer $C(x_{n-1}) = \langle x_1, x_3,\ldots, x_n\rangle$ is of codimension one.
This implies that $F(L) \subset C(x_{n-1})$ by [EO, Propositions 1.9 and 1.6].\\
Since $F(L)$ is not commutative we deduce that $L$ has no CP's, i.e. $\alpha(L) < c(L)$.  Hence $\alpha(L) = c(L) - 1 = q$.\\
From [O6, Theorem 51] we know that $Y(L) = k [x_n,f]$ where
$$f = 2x_1x_n + (-1)^{q+1} x_{q+1}^2 + 2 \sum\limits_{i=3}^q (-1)^i x_ix_{n-i+2}$$
We observe that $F(L) \not\subset H$ (as $x_1 \in F(L) \backslash H$).\\
Using Proposition 3.5 we may conclude that $M = S(H) Y(L) = k [x_{q+1}, \ldots, x_n, f]$ is a complete (it is even strongly complete) Poisson commutative subalgebra of $S(L)$, 
generated by elements of degree at most two.  Hence $L$ satisfies the Milovanov conjecture.\\
For the remaining cases it suffices to point out a CP.
\item[(3)] Suppose $L$ is of type $R_n$.\\
Basis of $L$ = $x_1,\ldots, x_n$, $n \geq 5$\\
Nonzero brackets: $[x_1, x_i] = x_{i+1}$, $i=2,\ldots, n-1$; $[x_2, x_j] = x_{j+2}$, $j=3,\ldots, n-2$.
One verifies that $i(L) = n-4$ and $c(L) = (n+n-4)/2 = n-2$.\\
On the other hand $C^2(L) = \langle x_3, x_4,\ldots, x_n\rangle$ is an abelian ideal of $L$ of dimension $n -2 = c(L)$ and so is a CP of $L$.  In particular, $\alpha(L) = n-2$.\\
The last equation can be obtained directly by using Proposition 8.4.  Indeed
$$\alpha(L) = n -\max \{i \mid [x_i, x_{i+1}]\neq 0\} = n-2$$
\item[(4)] Suppose $L$ is of type $W_n$.\\
Put $q = [(n+1)/2]$.  Then we claim that
$$H = C^{q-1}(L) = \langle x_q, x_{q+1},\ldots, x_n\rangle$$
is a CP of $L$.\\
First, we see that $\xi = x_n^\ast \in L^\ast$ is regular.\\
We distinguish 2 cases:
\begin{itemize}
\item[(4a)] $n$ is odd (i.e. $n = 2q-1$)\\
Then $L(\xi) = \langle x_n\rangle$ and so $i(L) = \dim L(\xi) =1$, $c(L) = q$.  Hence the claim follows from (1) of Proposition 8.5.  In particular, $\alpha(L) = q = [(n+2)/2]$
\item[(4b)] $n$ is even (i.e. $n = 2q$).\\
Then $L(\xi) = \langle x_q, x_n\rangle$ and so $i(L) = 2$, $c(L) = q + 1$.  By direct verification we see that
$$H = C^{q-1}(L) = \langle x_q, x_{q+1},\ldots, x_n\rangle$$
is abelian of dimension $2q - (q-1) = q+1 = c(L)$, i.e. $H$ is a CP of $L$ and $\alpha(L) = q+1 = [(n+2)/2]$.\hfill $\square$
\end{itemize}
\end{itemize}

{\bf Theorem 8.7.}\\
Let $L$ be an $n$-dimensional filiform Lie algebra.  Then $L$ satisfies the $R$-property if $n \leq 11$.\\
\ \\
{\bf Proof.} This is clear if $n \leq 7$ by (1) of Theorem 6.2 and if $n = 8$ by Corollary 6.3.\\
\ \\
{\bf (1) {\boldmath{$n = 9$}}}\\
The condition holds if $i(L) = 1$ by Proposition 8.5 and also if $i(L) = 7$ by Lemma 3.2.  So only the cases where $i(L) = 3$ or 5 remain.\\
By [B, Example 2.4.9] $L$ has an adapted basis $x_1,\ldots, x_9$ such that the brackets are given by (the undefined brackets are zero):\\
$m_{1i} = [x_1, x_i] = x_{i+1}$, $i = 2,\ldots, 8$.\\
$m_{23} = [x_2, x_3] = a_{25}x_5 + a_{26}x_6 + a_{27} x_7 + a_{28} x_8 + a_{29} x_9$\\
$m_{24} = [x_2, x_4] = a_{25}x_6 + a_{26}x_7 + a_{27} x_8 + a_{28} x_9$\\
$m_{25} = [x_2, x_5] = (a_{25}-a_{37})x_7 + (a_{26} - a_{38})x_8 + (a_{27} - a_{39}) x_9$\\
$m_{26} = [x_2, x_6] = (a_{25}-2a_{37})x_8 + (a_{26} - 2a_{38})x_9$\\
$m_{27} = [x_2, x_7] = (a_{25}-3a_{37} + a_{49}) x_9$\\
$m_{34} = [x_3, x_4] = a_{37}x_7 + a_{38}x_8 + a_{39} x_9$\\
$m_{35} = [x_3, x_5] = a_{37}x_8 + a_{38}x_9$\\
$m_{36} = [x_3, x_6] = (a_{37} - a_{49})x_9$\\
$m_{45} = [x_4, x_5] = a_{49}x_9$\\
The Jacobi identity holds if and only if the parameters $a_{ij} \in k$ satisfy the following equation:
$$a_{49} (2a_{25} + a_{37}) - 3a_{37}^2 = 0 \ \ \ (J)$$
The structure matrix $M = (m_{ij})$ with respect to the basis $x_1,\ldots, x_9$ is given in Table 4.
\newpage
\begin{center}
{\bf Table 4.}\\
\ \\
\begin{tabular}{c|ccccccccc}
      &$x_1$  &$x_2$     &$x_3$     &$x_4$     &$x_5$    &$x_6$    &$x_7$    &$x_8$ &$x_9$\\
\hline
$x_1$ &0      &$x_3$     &$x_4$     &$x_5$     &$x_6$    &$x_7$    &$x_8$    &$x_9$ &0\\
$x_2$ &$-x_3$ &0         &$m_{23}$  &$m_{24}$  &$m_{25}$ &$m_{26}$ &$m_{27}$ &0 &0\\
$x_3$ &$-x_4$ &$-m_{23}$ &0         &$m_{34}$  &$m_{35}$ &$m_{36}$ &0 &0 &0\\
$x_4$ &$-x_5$ &$-m_{24}$ &$-m_{34}$ &0         &$m_{45}$ &0        &0 &0 &0\\
$x_5$ &$-x_6$ &$-m_{25}$ &$-m_{35}$ &$-m_{45}$ &0        &0        &0 &0 &0\\
$x_6$ &$-x_7$ &$-m_{26}$ &$-m_{36}$ &0         &0        &0        &0 &0 &0\\
$x_7$ &$-x_8$ &$-m_{27}$ &0         &0         &0        &0        &0 &0 &0\\
$x_8$ &$-x_9$ &0         &0         &0         &0        &0        &0 &0 &0\\
$x_9$ &0      &0         &0         &0         &0        &0        &0 &0 &0
\end{tabular}
\end{center}

{\bf (1a) {\boldmath{$i(L) = 3$}}}\\
Then $\rank\ M = \dim L - i(L) = 9-3=6$ by formula (1), while $c(L) = (9+3)/2 = 6$. Clearly,
$$C^4(L) = \langle x_5, x_6, x_7, x_8, x_9\rangle$$
is an abelian ideal of $L$ of dimension $5 = c(L) - 1$, i.e. the R-property is valid.\\
\ \\
{\bf (1b) {\boldmath{$i(L) = 5$}}}\\
Then $\rank\ M = \dim L - i(L) = 9 -5 = 4$, while $c(L) = (9+5)/2 = 7$.  It suffices to show that $m_{45} = 0$, i.e. $a_{49} = 0$.  Because then
$$C^3(L) = \langle x_4, x_5, x_6, x_7, x_8, x_9\rangle$$
is an abelian ideal of $L$ of dimension $6 = c(L) -1$.\\
So, let us suppose that $a_{49} \neq 0$.  Now, consider the following submatrix çof $M$:
\begin{eqnarray*}
B = \left(\begin{array}{ccccc}
x_6 &x_7 &x_8 &x_9\\
m_{25} &m_{26} &m_{27} &0\\
m_{35} &m_{36} &0 &0\\
m_{45} &0 &0 &0\end{array}\right)
\end{eqnarray*}
Then $\rank\ B \leq 2$ since $\rank\ M = 4$ and by the special form of $M$.  This implies that $m_{26} = m_{36} = m_{27} = 0$, i.e.
$a_{25} = 2a_{37}$, $a_{37} = a_{49}$, $a_{25} - 3a_{37} + a_{49} = 0$.  Hence $a_{25} = 2a_{49}$, $a_{37} = a_{49}$.  Substitution in (J) gives 
$0 = a_{49} (4a_{49} + a_{49}) - 3a_{49}^2 = 2a_{49}^2$.  Contradiction.\\
\ \\
{\bf (2) {\boldmath{$n=10$}}}\\
The R-property holds if $i(L) = 2$ by Proposition 8.5 and also if $i(L) = 8$ by \\
Lemma 3.2.  So, the remaining cases are $i(L) = 4$ or 6.\\
By [B, p.70] $L$ has an adapted basis $x_1,\ldots, x_{10}$ such that the brackets are given by (the undefined brackets are zero):\\
$m_{1i} = [x_1, x_i] = x_{i+1}$, $i = 2,\ldots, 9$.\\
$m_{23} = [x_2, x_3] = a_{25}x_5 + a_{26}x_6 + a_{27} x_7 + a_{28} x_8 + a_{29} x_9 + a_{2,10} x_{10}$\\
$m_{24} = [x_2, x_4] = a_{25}x_6 + a_{26}x_7 + a_{27} x_8 + a_{28} x_9 + a_{29}x_{10}$\\
$m_{25} = [x_2, x_5] = (a_{25}-a_{37})x_7 + (a_{26} - a_{38})x_8 + (a_{27} - a_{39}) x_9 + (a_{28} - a_{3,10})x_{10}$\\
$m_{26} = [x_2, x_6] = (a_{25}-2a_{37})x_8 + (a_{26} - 2a_{38})x_9 + (a_{27} - 2a_{39}) x_{10}$\\
$m_{27} = [x_2, x_7] = (a_{25}-3a_{37} + a_{49}) x_9 + (a_{26} - 3a_{38} + \mu) x_{10}$\\
$m_{28} = [x_2, x_8] = (a_{25}-4a_{37} + 3a_{49}) x_{10}$\\
$m_{29} = [x_2,x_9] = - \lambda x_{10}$\\
$m_{34} = [x_3, x_4] = a_{37}x_7 + a_{38}x_8 + a_{39} x_9 + a_{3,10} x_{10}$\\
$m_{35} = [x_3, x_5] = a_{37}x_8 + a_{38}x_9 + a_{39} x_{10}$\\
$m_{36} = [x_3, x_6] = (a_{37} - a_{49})x_9 + (a_{38} - \mu) x_{10}$\\
$m_{37} = [x_3, x_7] = (a_{37} - 2a_{49}) x_{10}$\\
$m_{38} = [x_3, x_8] = \lambda x_{10}$\\
$m_{45} = [x_4, x_5] = a_{49}x_9 + \mu x_{10}$\\
$m_{46} = [x_4, x_6] = a_{49}x_{10}$\\
$m_{47} = [x_4, x_7] = -\lambda x_{10}$\\
$m_{5,6} = [x_5, x_6] = \lambda x_{10}$\\
The Jacobi identity holds if and only if the parameters satisfy the following equations:
\begin{eqnarray*}
\begin{array}{ll}
\lambda (2a_{25} - a_{37} - a_{49}) = 0 &{(J1)}\\
a_{49} (2a_{25} + a_{37}) - 3a_{37}^2 = 0 &{(J2)}\\
\lambda (2a_{27} + a_{39}) - \mu (2a_{25} + a_{37}) - 3a_{49} (a_{26} + a_{38}) + 7a_{37} a_{38} = 0 &{(J3)}
\end{array}
\end{eqnarray*}
The structure matrix $M = (m_{ij})$ with respect to the basis $x_1,\ldots, x_{10}$ is given in Table 5.
\newpage
\begin{center}
{\bf Table 5.}\\
\ \\
\begin{tabular}{c|cccccccccc}
      &$x_1$  &$x_2$     &$x_3$     &$x_4$     &$x_5$    &$x_6$    &$x_7$    &$x_8$    &$x_9$             &$x_{10}$\\
\hline
$x_1$ &0      &$x_3$     &$x_4$     &$x_5$     &$x_6$    &$x_7$    &$x_8$    &$x_9$    &$x_{10}$          &0\\
$x_2$ &$-x_3$ &0         &$m_{23}$  &$m_{24}$  &$m_{25}$ &$m_{26}$ &$m_{27}$ &$m_{28}$ &$-\lambda x_{10}$ &0\\
$x_3$ &$-x_4$ &$-m_{23}$ &0         &$m_{34}$  &$m_{35}$ &$m_{36}$ &$m_{37}$ &$\lambda x_{10}$ &0         &0\\
$x_4$ &$-x_5$ &$-m_{24}$ &$-m_{34}$ &0         &$m_{45}$ &$m_{46}$ &$-\lambda x_{10}$        &0 &0 &0\\
$x_5$ &$-x_6$ &$-m_{25}$ &$-m_{35}$ &$-m_{45}$ &0        &$\lambda x_{10}$        &0 &0 &0 &0\\
$x_6$ &$-x_7$ &$-m_{26}$ &$-m_{36}$ &$-m_{46}$ &$-\lambda x_{10}$        &0 &0        &0 &0 &0\\
$x_7$ &$-x_8$ &$-m_{27}$ &$-m_{37}$ &$\lambda x_{10}$         &0         &0        &0        &0 &0 &0\\
$x_8$ &$-x_9$ &$-m_{28}$ &$-\lambda x_{10}$         &0         &0         &0        &0        &0 &0 &0\\
$x_9$ &$-x_{10}$ &$\lambda x_{10}$    &0         &0         &0         &0        &0        &0 &0 &0\\
$x_{10}$ &0 &0 &0 &0 &0 &0 &0 &0 &0 &0
\end{tabular}
\end{center}
{\bf (2a) {\boldmath{$i(L) =4$}}}\\
Then $\rank\ M = \dim L - i(L) = 10 - 4 = 6$ by formula (1), which implies that $\lambda = 0$.  Also, $c(L) = (10+4)/2 = 7$.\\
Clearly,
$$C^4(L) = \langle x_5, x_6, x_7, x_8, x_9, x_{10}\rangle$$
is an abelian ideal of $L$ of dimension $6 = c(L) -1$, i.e. the R-property holds for $L$.\\
\ \\
{\bf (2b) {\boldmath{$i(L) = 6$}}}\\
Then $\rank\ M = 10-6 = 4$.  Hence $\lambda = 0$.  Also, $c(L) = (10+6)/2 = 8$.  It suffices to show that $m_{45} = [x_4, x_5] = 0$.  Because then $m_{46} = [x_4,x_6] = [x_1, m_{45}] = 0$ and thus
$$C^3(L) = \langle x_4, x_5, x_6, x_7, x_8, x_9, x_{10}\rangle$$
is an abelian ideal of $L$ of dimension $7 = c(L) - 1$.  Hence the R-property is satisfied.\\
So, let us suppose that $m_{45} \neq 0$.\\
Because $\rank\ M = 4$ and by the special form of $M$ we see that
\begin{eqnarray*}
0 = \det \left(\begin{array}{ccc}
x_6 &x_7 &x_{10}\\
m_{35} &m_{36} &0\\\
m_{45} &m_{46} &0\end{array}\right) = (m_{35}m_{46} - m_{36}m_{45}) x_{10}
\end{eqnarray*}
Hence $m_{35}m_{46} = m_{36}m_{45}$, i.e.
$$(a_{37}x_8 + a_{38}x_9 + a_{39}x_{10}) a_{49}x_{10} = [(a_{37} - a_{49})x_9 + (a_{38} - \mu) x_{10}] (a_{49}x_9 + \mu x_{10})$$
From the identification of the coefficients we obtain:
$$a_{37}a_{49} = 0,\ \ a_{37}a_{49} - a_{49}^2 = 0,\ \ a_{39}a_{49} = (a_{38} - \mu)\mu$$
It follows that $a_{49} = 0$.  Hence $\mu \neq 0$ (since $0 \neq m_{45} = \mu x_{10}$) and so $a_{38} = \mu$.  From (J2) we get $a_{37} = 0$, indeed $3a_{37}^2 = a_{49} (2a_{25} + a_{37}) = 0$.  
Moreover, (J3) now reduces to $\mu(2a_{25} + a_{37}) = 0$ and thus $a_{25} = 0$.\\
Next, $\rank\ M = 4$ and $\lambda = 0$ imply that
\begin{eqnarray*}
-m_{45}m_{27}x_{10} = \det \left(\begin{array}{ccc}
x_6 &x_8 &x_{10}\\
m_{25} &m_{27} &0\\
m_{45} &0 & 0\end{array}\right) = 0
\end{eqnarray*}
Therefore $m_{27} = 0$.  In particular $a_{26} - 3a_{38} + \mu = 0$, i.e. $a_{26} = 2\mu$.\\
Since $\rank\ M = 4$, the following $6 \times 6$ submatrix of $M$ has a zero determinant.
\begin{eqnarray*}
\left(\begin{array}{cccccc}
0 &x_3 &x_4 &x_5 &x_6 &x_{10}\\
-x_3 &0 &m_{23} &m_{24} &m_{25} &0\\
-x_4 &-m_{23} &0 &m_{34} &m_{35} &0\\
-x_5 &-m_{24} &-m_{34} &0 &m_{45} &0\\
-x_6 &-m_{25} &-m_{35} &-m_{45} &0 &0\\
-x_{10} &0 &0 &0 &0 &0\end{array}\right)
\end{eqnarray*}
It follows that
\begin{eqnarray*}
0 = \det \left(\begin{array}{cccc}
0 &m_{23} &m_{24} &m_{25}\\
-m_{23} &0 &m_{34} &m_{35}\\
-m_{24} &-m_{34} &0 &m_{45}\\
-m_{25} &-m_{25} &-m_{45} &0\end{array}\right) = (m_{24} m_{35} - m_{25}m_{34} - m_{23} m_{45})^2
\end{eqnarray*}
Therefore, $m_{24} m_{35} = m_{25} m_{34} + m_{23} m_{45}$.\\
Taking into account that $a_{49} = a_{37} = a_{25} = 0$, $a_{38} = \mu$, $a_{26} = 2\mu$. we obtain:
\\
$(2\mu x_7 + a_{27}x_8 + a_{28}x_9 + a_{29} x_{10}) (\mu x_9 + a_{39}x_{10}) =$\\
$[\mu x_8 + (a_{27} - a_{39}) x_9 + (a_{28} - a_{3,10}) x_{10}](\mu x_8 + a_{39} x_9 + a_{3,10} x_{10}) +$\\
$(2\mu x_6 + a_{27} x_7 + a_{28} x_8 + a_{29} x_9 + a_{2,10}x_{10} ) \mu x_{10}$\\
Comparing the coefficients of $x_7x_9$ (or of $x_6 x_{10})$ of both sides, we get $2\mu^2 = 0$.  Contradiction.\\
\ \\
{\bf (3) {\boldmath{$n=11$}}}\\
The proof, which uses the same approach as above, is quite long and therefore it will be omitted. \hfill $\square$\\
\newpage
{\bf {\large 9. The Poisson center and Milovanov's conjecture for filiform Lie algebras of dimension {\boldmath{$\leq 8 (k = \Bbb C)$}}}}\\
Let $\mathfrak{g}$ be filiform of dimension $n \geq 8$.  By Theorem 8.7 it satisfies the R-property.  The Poisson center $Y(\mathfrak{g})$ of $\mathfrak{g}$ is not only an interesting object
in its own right, but for us it is a sueful tool in the construction, as outlined in 1.3, of a complete Poisson commutative subalgebra $M$ of $S(\mathfrak{g})$.  It turns out that $M$ (or sometimes a slight enlargement of $M$)
will be generated by elements of degree at most two, i.e. the Milovanov conjecture is valid for $\mathfrak{g}$.
If $n \leq 7$ this has been established already [O4, O5].  Assume $n = 8$. We recall from [O6, Proposition 50] that:
\begin{center}
$\mathfrak{g}$ is coregular $\Leftrightarrow i(\mathfrak{g}) = 2$
\end{center}
In case $i(\mathfrak{g}) = 2$ we will exhibit algebraically independent generators of $Y(\mathfrak{g})$.\\
However, if $i(\mathfrak{g}) \geq 4$, i.e. $Y(\mathfrak{g})$ is not polynomial, exhibiting the generators of $Y(\mathfrak{g})$ becomes quite complicated (see e.g. [O5, Example 27]).  Therefore we will only provide
algebraically independent generators of the quotient field $Q(Y(\mathfrak{g}))$ of $Y(\mathfrak{g})$ by using a technique due to Dixmier [O5, Theorem 31].\\
\ \\
Let us now illustrate this by the following:\\
{\bf Example 9.1.} $\mathfrak{g} = \mathfrak{g}_{8,9}(\lambda)$\\
{\bf (i) {\boldmath{$\lambda \neq 1$}}.}\\
Basis: $x_1,x_2,\ldots, x_8$\\
Nonzero brackets: $[x_1,x_j] = x_{j+1}$, $j = 2,\ldots, 7$, $[x_2, x_3] = \lambda x_6 + x_7$, $[x_2, x_4] = \lambda x_7 + x_8$, $[x_2, x_5] = (\lambda-1) x_8$, $[x_3, x_4] = x_8$.\\
Clearly $\xi = x_8^\ast \in \mathfrak{g}^\ast$ is regular, $\mathfrak{g}(\xi) = \langle x_6, x_7\rangle$ and $F(\mathfrak{g}) = \langle x_6, x_7, x_8\rangle$.  Therefore, $i(\mathfrak{g}) = \dim \mathfrak{g}(\xi) = 2$ and $c(\mathfrak{g}) = (8+2)/2 = 5$.  Since $x_1, \ldots, x_8$ is an adapted basis
of $\mathfrak{g}$ we obtain from Proposition 8.4 that
$$\alpha(\mathfrak{g}) = 8 - \max \{j \mid [x_j, x_{j+1}] \neq 0\} = 8-3=5$$
and $\mathfrak{h} = C^3 (\mathfrak{g}) = \langle x_4, x_5, x_6, x_7, x_8\rangle$ is abelian with $\dim \mathfrak{h} = 5 = c(\mathfrak{g})$.  Consequently, $\mathfrak{h}$ is a CP of $\mathfrak{g}$.\\
One verifies that $p_{\mathfrak{g}} = x_8^2$ and $x_8$, $f = 2x_6x_8 - x_7^2 \in Y(\mathfrak{g})$.\\
Using the claim of the proof of [O6, Theorem 45] we get at once that $Y(\mathfrak{g}) = k[x_8,f]$.  But this can also be seen as an application of [JS, 5.7], [O6, Theorem 29], since $i(\mathfrak{g}) = 2$ and
$$\deg x_8 + \deg f = 3 = 5-2 = c(\mathfrak{g}) - \deg p_{\mathfrak{g}}$$
Finally, $M = k [x_4, x_5, x_6, x_7, x_8]$ is a strongly complete Poisson commutative subalgebra of $S(\mathfrak{g})$ because $\mathfrak{h} = \langle x_4, x_5, x_6, x_7, x_8\rangle$
 is a CP of $\mathfrak{g}$. So Milovanov's conjecture is trivially satisfied.\\
\ \\
{\bf (ii) {\boldmath{$\lambda = 1$}}.}\\
Basis: $x_1,x_2,\ldots, x_8$\\
Nonzero brackets: $[x_1,x_j] = x_{j+1}$, $j = 2,\ldots, 7$, \\
$[x_2, x_3] = x_6 + x_7$, $[x_2, x_4] =  x_7 + x_8$, $[x_3, x_4] = x_8$.\\
Clearly $\xi = x_8^\ast$ is regular, $\mathfrak{g}(\xi) = \langle x_2 -x_3, x_5, x_6, x_8\rangle$ and $F(\mathfrak{g}) = \langle x_2, x_3, x_4, x_5, x_6, x_7, x_8\rangle$, 
which is not abelian, so there are no CP's.\\
$i(\mathfrak{g}) = \dim \mathfrak{g}(\xi) = 4 > 2$, so $Y(\mathfrak{g})$ is not polynomial, $c(\mathfrak{g}) = 6$, $p_{\mathfrak{g}} = 1$.\\
As above we see that $\alpha(\mathfrak{g}) = 5 = c(\mathfrak{g})-1$ and also that
$$\mathfrak{h} = C^3 (\mathfrak{g}) = \langle x_4, x_5, x_6, x_7, x_8\rangle$$
is abelian of dimension 5.  The following invariants:\\
$x_8$,\\
$f_1 = 2x_6x_8 - x_7^2$,\\
$f_2 = 3x_5x_8^2 - 3x_6 x_7x_8 + x_7^3$,\\
$f_3 = 10(x_2 - x_3) x_8^2 - 10x_3 x_7 x_8 + 10x_4 x_6 x_8 + 10x_4 x_7 x_8 - 5x_5^2x_8 - 2x_5 x_6x_8 + 2x_6^2 x_7 - 4x_5 x_7^2$\\
$= 10 (x_2 x_8 - x_3x_8 - x_3x_7) x_8 + g$,\\
where $g = 10x_4x_6x_8 + 10x_4x_7x_8 - 5x_5^2x_8 - 2x_5x_6x_8 + 2x_6^2x_7 - 4x_5 x_7^2 \in S(\mathfrak{h})$,\\
are algebraically independent generators of $Q(Y(\mathfrak{g}))$ by [O5, Theorem 31].\\
Note that $f_3 \in Y(\mathfrak{g})\backslash S(\mathfrak{h})$.  By (i) of Remark 3.4
\begin{eqnarray*}
S(\mathfrak{h}) k[f_3] &=& k[x_4,x_5,x_6,x_7, x_8, f_3]\\
&=& k[x_4, x_5, x_6, x_7, x_8, (x_2x_8 - x_3x_8 - x_3x_7)x_8]
\end{eqnarray*}
is a complete, Poisson commutative subalgebra of $S(\mathfrak{g})$.\\
The same is true for
$$M = k[x_4,x_5,x_6,x_7,x_8, x_2x_8 - x_3x_8 - x_3x_7]$$
which satisfies the conditions of Milovanov's conjecture.  $M$ is also strongly complete.  Indeed, consider the Jacobian locus $J$ of the generators of $M$. \\ 
Then $J = \{\xi \in \mathfrak{g}^\ast \mid \xi (x_7) = 0 = \xi (x_8)\}$.\\
Clearly $\codim\ J = 2$.  The result then follows from [PPY, Theorem 1.1] combined with [PY, 2.1].\\
\newpage
{\bf 9.2. The list of filiform Lie algebras of dimension {\boldmath{$n \leq 8$}}}\\
Our primary aim is to exhibit generators of the Poisson center $Y(\mathfrak{g})$ (or in some cases only of its quotient field $Q(Y(\mathfrak{g}))$) of each member $\mathfrak{g}$ if the list\\
Secondly we will produce, among other things, a polynomial, complete Poisson commutative subalgebra $M$ of $S(\mathfrak{g})$, generated by elements of degree at most two.
\begin{itemize}
\item[$\bullet$] If $n \leq 7$ we will simply select the filiform Lie algebras from the list of all indecomposable nilpotent Lie algebras of dimension at most seven [O4, O5]. See also [GK, pp.58-62].
\item[$\bullet$] If $n = 8$ our list is based on the classification  of [GJK].  See also [AG].
\end{itemize}

{\bf Notation and abbreviations:}\\
$x_1,\ldots, x_n$ will be a basis of $\mathfrak{g}$, $n \leq 8$.\\
SQ.I. = square integrable, $i=i(\mathfrak{g})$, $c=c(\mathfrak{g})$, $p=p_{\mathfrak{g}}$, $F = F(\mathfrak{g})$, $\alpha = \alpha(\mathfrak{g})$, $\mathfrak{h}$ is an abelian ideal of $\mathfrak{g}$
with $\dim \mathfrak{h} = \alpha(\mathfrak{g})$, $C^j = C^j(\mathfrak{g})$, $Y = Y(\mathfrak{g})$, $Q(Y) = Q(Y(\mathfrak{g}))$.\\
\ \\
{\bf (i) {\boldmath{$n \leq 5$}}}
\begin{itemize}
\item[1.] $\mathfrak{g}_3$ (=1 of [O4]) $=L_3$\\
$[x_1, x_2] = x_3$.\\
SQ.I.\quad $i = 1$, $c=2$, $p= x_3$, $F = \langle x_3\rangle$, $\alpha = 2$, $\mathfrak{h} = \langle x_2, x_3\rangle = CP$,\\
$Y = k[x_3]$, $M = k[x_2, x_3]$.

\item[2.] $\mathfrak{g}_4$ (= 2 of [O4]) $= L_4$\\
$[x_1, x_2] = x_3$, $[x_1, x_3] = x_4$.\\
$i = 2$, $c=3$, $p=1$, $\alpha = 3$, $F = \langle x_2, x_3, x_4) = \mathfrak{h} = CP$,\\
$Y = k[x_4, x_3^2 - 2x_2x_4]$, $M = k[x_2, x_3, x_4]$.

\item[3.] $\mathfrak{g}_{5,6}$ (=6 of [O4]) $= R_5$\\
$[x_1,x_2] = x_3$, $[x_1, x_3] = x_4$, $[x_1, x_4] = x_5$, $[x_2, x_3] = x_5$\\
SQ.I., $i=1$, $c=3$, $p = x_5^2$, $F = \langle x_5\rangle$, $\alpha = 3$, $\mathfrak{h} = C^2 =\langle x_3, x_4, x_5\rangle= CP$,\\
$Y = k[x_5]$, $M = k[x_3, x_4, x_5]$.

\item[4.] $\mathfrak{g}_{5,5}$ (=8 of [O4]) $=L_5$, not coregular\\
$[x_1, x_2] = x_3$, $[x_1, x_3] = x_4$, $[x_1, x_4] = x_5$.\\
$i = 3$, $c = 4$, $p=1$, $\alpha = 4$, $F = \langle x_2,x_3,x_4,x_5\rangle = \mathfrak{h} = CP$\\
$Y = k[x_5, f_1, f_2, f_3]$, $f_1 = 2x_3x_5 - x_4^2$, $f_2 = 3x_2x_5^2 - 3x_3 x_4 x_5 + x_4^3$ \\
$f_3 = 9x_2^2x_5^2 - 18x_2x_3x_4x_5 + 6x_2x_4^3 + 8x_3^3x_5 - 3x_3^2x_4^2$\\
Relation: \quad $f_1^3 + f_2^2 - x_5^2f_3 = 0$, $Q(Y) = k(x_5, f_1, f_2)$,\\
$M = k[x_2, x_3, x_4, x_5]$.
\end{itemize}

{\bf (ii) {\boldmath{$n = 6$}}}
\begin{itemize}
\item[5.] $\mathfrak{g}_{6,18}$ (=21 of O4]) $= Q_6$\\
$[x_1, x_2] = x_3$, $[x_1, x_3] = x_4$, $[x_1, x_4] = x_5$, $[x_2, x_5] = x_6$, $[x_3, x_4] = -x_6$.\\
$i = 2$, $c = 4$, $p = x_6$, $F =\langle x_1,x_3,x_4,x_5,x_6\rangle$; no CP's, $\alpha = 3$,\\
$\mathfrak{h} = C^3 = \langle x_4, x_5, x_6\rangle$, $Y = k[x_6, x_4^2 - 2x_3x_5 - 2x_1x_6]$.\\
$M = k[x_4, x_5, x_6, x_3x_5 + x_1x_6]$.

\item[6.] $\mathfrak{g}_{6,17}$ (=26 of [O4])\\
$[x_1, x_2] = x_3$, $[x_1, x_3] = x_4$, $[x_1, x_4] = x_5$, $[x_1, x_5] = x_6$, $[x_2,x_3] = x_6$.\\
$i = 2$, $c = 4$, $p = x_6$, $F = \langle x_4,x_5,x_6\rangle$, $\alpha = 4$, $\mathfrak{h} = C^2 = \langle x_3, x_4, x_5, x_6\rangle = CP$,
$Y = k[x_6, x_5^2 - 2x_4x_6]$, $M = k[x_3,x_4,x_5,x_6].$

\item[7.] $\mathfrak{g}_{6,19}$ (=27 of [O4]) $= R_6$\\
$[x_1, x_2] = x_3$, $[x_1, x_3] = x_4$, $[x_1, x_4] = x_5$, $[x_1, x_5] = x_6$, $[x_2, x_3] = x_5$,\\
$[x_2, x_4] = x_6$.\\
$i = 2$, $c = 4$, $p = 1$, $\alpha = 4$, $F = \langle x_3, x_4, x_5,x_6\rangle = \mathfrak{h} = C^2 = CP$,\\
$Y = k[x_6, x_5^3 -3x_4x_5x_6 + 3x_3x_6^2]$, $M = k[x_3,x_4,x_5,x_6]$.

\item[8.] $\mathfrak{g}_{6,20}$ (=28 of [O4])\\
$[x_1, x_2] = x_3$, $[x_1, x_3] = x_4$, $[x_1, x_4] = x_5$, $[x_2, x_3] = x_5$, $[x_2, x_5] = x_6$,\\
$[x_3, x_4] = -x_6$.\\
$i = 2$, $c = 4$, $p = 1$, $F = \langle x_1, x_3, x_4, x_5,x_6\rangle$, no CP's, $\alpha = 3$,\\
$\mathfrak{h} = C^3 = \langle x_4, x_5, x_6\rangle$,\\
$Y = k[x_6, 2x_5^3 + 3x_4^2 x_6 -6x_3x_5x_6 - 6x_1x_6^2]$, $M = k[x_4, x_5, x_6, x_3x_5 + x_1x_6]$.

\item[9.] $\mathfrak{g}_{6,16}$ (=25 of [O4]) $= L_6$, not coregular\\
$[x_1, x_2] = x_3$, $[x_1, x_3] = x_4$, $[x_1, x_4] = x_5$, $[x_1, x_5] = x_6$.\\
$i = 4$, $c = 5$, $p = 1$, $F = \langle x_2, x_3,x_4, x_5,x_6\rangle = CP= \mathfrak{h}$, $\alpha = 5$,\\
$Y = k[x_6, f_1, f_2, f_3, f_4]$, $f_1 = x_5^2 - 2x_4x_6$, $f_2 = x_5^3 - 3x_4x_5x_6 + 3x_3x_6^2,\\
f_3 = x_4^2 + 2x_2x_6 - 2x_3x_5,\\
f_4 = 2x_4^3 + 6x_2x_5^2 + 9x_3^2x_6 - 12x_2x_4x_6 - 6x_3x_4x_5$.\\
Relation: \quad $f_1^3 - f_2^2 - 3x_6^2f_1f_3 + x_6^3f_4 = 0$, $Q(Y) = k(x_6, f_1, f_2,f_3)$.\\
$M = k[x_2, x_3, x_4, x_5, x_6]$.
\end{itemize}

{\bf {\boldmath{$n = 7$}}}\\
{\bf 7.a {$\boldmath{\mathfrak{g}}$} is coregular}
\begin{itemize}
\item[10.] $\mathfrak{g}_{7,1.1(i_\lambda)}$, $\lambda \neq 0,1$ (= 30 of [O4])\\
$[x_1, x_2] = x_3$, $[x_1, x_3] = x_4$, $[x_1, x_4] = x_5$, $[x_1, x_5] = x_6$, $[x_1, x_6] = x_7$,\\
$[x_2, x_3] = x_5$, $[x_2, x_4] = x_6$, $[x_2, x_5] = \lambda x_7$, $[x_3, x_4] = (1-\lambda) x_7$.\\
SQ.I., $i = 1$, $c = 4$, $p = x_7^3$, $F = \langle x_7\rangle$, $\alpha = 4$, $\mathfrak{h} = C^3 = \langle x_4, x_5, x_6, x_7\rangle = CP$,\\
$Y = k[x_7]$, $M = k[x_4, x_5, x_6, x_7]$.

\item[11.] $\mathfrak{g}_{7,1.1(ii)}$ (= 31 of [O4])\\
$[x_1, x_2] = x_3$, $[x_1, x_3] = x_4$, $[x_1, x_4] = x_5$, $[x_1, x_5] = x_6$, $[x_1, x_6] = x_7$\\
$[x_2, x_5] = x_7$, $[x_3, x_4] = -x_7$.\\
SQ.I., $i = 1$, $c = 4$, $p = x_7^3$, $F = \langle x_7\rangle$, $\alpha = 4$, $\mathfrak{h} = C^3 = \langle x_4, x_5, x_6, x_7\rangle = CP$,\\
$Y = k[x_7]$, $M = k[x_4, x_5, x_6, x_7]$.

\item[12.] $\mathfrak{g}_{7,0.1}$ (= 83 of [O5])\\
$[x_1,x_2] = x_3$, $[x_1, x_3] = x_4$, $[x_1, x_4] = x_5$, $[x_1,x_5] = x_6$, $[x_1, x_6] = x_7$,\\
$[x_2, x_3] = x_6$, $[x_2, x_4] = x_7$, $[x_2, x_5] = x_7$, $[x_3, x_4] = -x_7$.\\
SQ.I., $i = 1$, $c=4$, $p = x_7^3$, $F = \langle x_7\rangle$, $\alpha = 4$, $\mathfrak{h} = C^3 = \langle x_4, x_5, x_6, x_7\rangle$, $= CP$,\\
$Y = k[x_7]$, $M = k[x_4, x_5, x_6, x_7]$.

\item[13.] $\mathfrak{g}_{7,1.4}$ (= 106 of [O5])\\
$[x_1,x_2] = x_3$, $[x_1, x_3] = x_4$, $[x_1, x_4] = x_5$, $[x_1,x_5] = x_6$, $[x_1, x_6] = x_7$,\\
$[x_2, x_3] = x_6$, $[x_2, x_4] = x_7$.\\
$i = 3$, $c=5$, $p=1$, $\alpha = 5$, $F = \langle x_3, x_4, x_5, x_6, x_7\rangle = C^2 = \mathfrak{h} = CP$,\\
$Y = k[x_7, x_5^2 - 2x_4x_6 + 2x_3x_7, x_6^2 - 2x_5x_7]$, $M = k[x_3, x_4, x_5, x_6, x_7]$.
\end{itemize}

{\bf 7.b {$\boldmath{\mathfrak{g}}$} is not coregular}
\begin{itemize}

\item[14.] $\mathfrak{g}_{7,0.2}$ (= 153 of [O5])\\
$[x_1,x_2] = x_3$, $[x_1, x_3] = x_4$, $[x_1, x_4] = x_5$, $[x_1,x_5] = x_6$, $[x_1, x_6] = x_7,$ \\
$[x_2, x_3] = x_5+x_7, [x_2, x_4] = x_6, [x_2, x_5] = x_7$.\\
$i = 3$, $c= 5$, $p = 1$, $\alpha = 5$, $F = \langle x_3, x_4, x_5, x_6, x_7\rangle = \mathfrak{h} = C^2 = CP$,\\
$Y = k[x_7, f, g, h]$, $f = x_6^3 - 3x_5x_6x_7 + 3x_4x_7^2$,\\
$g = x_6^4 - 4x_5x_6^2x_7 + 2x_5^2x_7^2 + 4x_4x_6x_7^2 - 2x_6^2x_7^2 - 4x_3x_7^3 + 4x_5x_7^3$,\\
$h = (f^4-g^3-6x_7^2f^2g) / x_7^3$, relation: $f^4 - g^3 - 6x_7^2f^2g - x_7^3 h = 0$,\\
$Q(Y) = k(x_7, f, g)$, $M = k[x_3, x_4, x_5, x_6, x_7]$.

\item[15.] $\mathfrak{g}_{7,0.3}$ (= 141 of [O5])\\
$[x_1,x_2] = x_3$, $[x_1, x_3] = x_4$, $[x_1, x_4] = x_5$, $[x_1,x_5] = x_6$, $[x_1, x_6] = x_7,$\\
$[x_2, x_3] = x_6+x_7, [x_2, x_4] = x_7$.\\
$i = 3$, $c=5$, $p = 1$, $\alpha = 5$, $F = \langle x_3, x_4, x_5, x_6, x_7\rangle = \mathfrak{h} = C^2 = CP$, \\
$Y = k[x_7, f, g, h]$, $f = x_6^2 - 2x_5x_7$,\\
$g = 2x_6^3 - 3x_5^2x_7 + 6x_4x_6x_7 -6x_5x_6x_7 - 6x_3x_7^2 + 6x_4x_7^2$,\\
$h = (4f^3-g^2) / x_7$, \\
relation: $4f^3 - g^2 - x_7h = 0$, $Q(Y) = k(x_7, f, g)$, $M = k[x_3, x_4, x_5, x_6, x_7]$.

\item[16.] $\mathfrak{g}_{7,1.6}$ (= 137 of [O5])\\
$[x_1,x_2] = x_3$, $[x_1, x_3] = x_4$, $[x_1, x_4] = x_5$, $[x_1,x_5] = x_6$, $[x_1, x_6] = x_7$,\\
$[x_2,x_3] = x_7$.\\
$i = 3$, $c=5$, $p = x_7$, $F = \langle x_4, x_5, x_6, x_7\rangle$, $\alpha = 5$, $\mathfrak{h} = C^2 = \langle x_3, x_4, x_5, x_6, x_7\rangle$ = CP,\\
$Y = k[x_7, f, g, h]$, $f = x_6^2 - 2x_5x_7$, $g = x_6^3 - 3x_5x_6x_7 + 3x_4x_7^2$,\\
$h = (f^3-g^2) / x_7^2$, relation: $f^3 - g^2 - x_7^2h = 0$,\\
$Q(Y) = k(x_7, f, g)$, $M = k[x_3, x_4, x_5, x_6, x_7]$.

\item[17.] $\mathfrak{g}_{7,1.1(i_\lambda),\lambda =1)}$ (= 151 of [O5]) $= R_7$\\
$[x_1,x_2] = x_3$, $[x_1, x_3] = x_4$, $[x_1, x_4] = x_5$, $[x_1,x_5] = x_6$, $[x_1, x_6] = x_7,$\\
$[x_2, x_3] = x_5, [x_2, x_4] = x_6, [x_2,x_5] = x_7$.\\
$i = 3$, $c=5$, $p = 1$, $\alpha = 5$, $F = \langle x_3, x_4, x_5, x_6, x_7\rangle = \mathfrak{h} = C^2 = CP$, \\
$Y = k[x_7, f, g, h]$, $f = x_6^3 - 3x_5x_6x_7 + 3x_4x_7^2$,\\
$g = x_6^4 - 4x_5x_6^2x_7 + 2x_5^2x_7^2 + 4x_4x_6x_7^2 - 4x_3x_7^3$, $h = (f^4-g^3) / x_7^3$, \\
relation: $f^4 - g^3 - x_7^3h = 0$, $Q(Y) = k(x_7, f, g)$,\\
$M = k[x_3, x_4, x_5, x_6, x_7]$.

\item[18.] $\mathfrak{g}_{7,1.1(i_\lambda)}, \lambda = 0$ (= 155 of [O5])\\
$[x_1,x_2] = x_3$, $[x_1, x_3] = x_4$, $[x_1, x_4] = x_5$, $[x_1,x_5] = x_6$, $[x_1, x_6] = x_7,$\\
$[x_2, x_3] = x_5, [x_2, x_4] = x_6, [x_3, x_4] = x_7$.\\
$i = 3$, $c=5$, $p=1$, $F = \langle x_2, x_3, x_4, x_5, x_6, x_7\rangle$, no $CP$'s,\\
$\alpha = 4$, $\mathfrak{h} =\langle x_4, x_5, x_6, x_7\rangle = C^3$\\
$Y = k[x_7, f, g, h]$, $f = x_6^2 - 2x_5x_7$,\\
$g = 2x_6^5 - 10x_5x_6^3 x_7 + 15x_5^2x_6x_7^2 - 15x_4x_5x_7^3 + 15x_3x_6x_7^3 - 15x_2x_7^4$,\\
$h = (4f^5-g^2) / x_7^3$, relation: $4f^5 - g^2 - x_7^3h = 0$,\\
$Q(Y) = k(x_7, f, g)$, \\
$M = k[x_4, x_5, x_6, x_7, x_3x_6 - x_2x_7]$.

\item[19.] $\mathfrak{g}_{7,2.3}$ (= 159 of [O5]) $= L_7$\\
$[x_1,x_2] = x_3$, $[x_1, x_3] = x_4$, $[x_1, x_4] = x_5, [x_1, x_5] = x_6, [x_1, x_6] = x_7$.\\
$i = 5$, $c=6$, $p=1$, $\alpha = 6$, $F = \langle x_2, x_3, x_4, x_5, x_6, x_7\rangle = \mathfrak{h} = CP$, \\
$Y = k[f_1, f_2, \ldots, f_{23}]$, $f_1 = x_7$, $f_2 = x_6^2 - 2x_5x_7$,\\
$f_3 = x_6^3 - 3x_5x_6x_7 + 3x_4x_7^2$, $f_4 = x_5^2 - 2x_4x_6 + 2x_3x_7$,\\
$f_5 = 2x_4x_6^2 - x_5^2x_6 + x_4x_5x_7 - 5x_3x_6x_7 + 5x_2x_7^2$ are algebraically independent over $k$, 
$Q(Y) = k(f_1, f_2, f_3, f_4, f_5)$, $M = k[x_2, x_3, x_4, x_5, x_6, x_7]$.
\end{itemize}

{\bf {\boldmath{$n = 8$}}}\\
{\bf 8.a {$\boldmath{\mathfrak{g}}$} is coregular (i.e. {$\boldmath{i(g) = 2}$}})
\begin{itemize}
\item[20.] $\mathfrak{g}_{8.1}(\lambda)$\\
$[x_1, x_j] = x_{j+1}$, $j = 2,\ldots, 7$,\\
$[x_2, x_3] = x_5 + \lambda x_6$, $[x_2,x_4] = x_6 + \lambda x_7$, $[x_2, x_5] = 3x_7 + \lambda x_8$, $[x_2, x_6] = 5x_8$,\\
$[x_2, x_7] = x_8$, $[x_3, x_4] = -2x_7$, $[x_3, x_5] = -2x_8$, $[x_3, x_6] = -x_8$, $[x_4, x_5] = x_8$.\\
$i=2$, $c=5$, $p=1$, $F = \langle x_1 - x_2, x_3, x_4, x_5, x_6, x_7, x_8\rangle$, no CP's,\\
$\alpha = 4$, $\mathfrak{h} = C^4 = \langle x_5, x_6, x_7, x_8\rangle$, $Y = k[x_8, f]$,\\
$f = x_5^2 x_8^2 - 2(\lambda -10) x_5 x_7x_8^2 + (\lambda + 2) x_6^2 x_8^2 - 12x_6 x_7^2 x_8 - 4x_5 x_7^2 x_8$\\
$+ 2x_6^2 x_7 x_8 - 2x_5 x_6 x_8^2 + 3x_7^4 + 2x_8^2 g$, where\\
$g = x_3x_7 + 5x_4x_7 - x_4x_6 + (\lambda - 10) x_4x_8 + x_1x_8 - x_2x_8 - 5x_3x_8$,\\
$M = k[x_5, x_6, x_7, x_8, g]$.

\item[21.] $\mathfrak{g}_{8.2}(\lambda \neq 0)$\\
$[x_1, x_j] = x_{j+1}$, $j = 2,\ldots, 7$,\\
$[x_2, x_3] = x_6 + \lambda x_7$, $[x_2,x_4] = x_7 + \lambda x_8$, $[x_2, x_5] = x_8$, $[x_2, x_7] = x_8$,\\
$[x_3, x_6] = -x_8$, $[x_4, x_5] = x_8$.\\
$i=2$, $c=5$, $p=x_8$, $F = \langle x_1 - x_2, x_3, x_4, x_5, x_6, x_7, x_8\rangle$, no CP's,\\
$\alpha = 4$, $\mathfrak{h} = C^4 = \langle x_5, x_6, x_7, x_8\rangle$, $Y = k[x_8, f]$,\\
$f = 3x_5^2 x_8 + 3x_6^2x_8 - 2\lambda x_7^3 - 6x_5x_7x_8 + 6\lambda x_6 x_7 x_8 - 6\lambda x_5 x_8^2 + 6x_8 g$\\
where $g = (x_1-x_2)x_8 - x_4x_6 + x_3x_7 + x_4x_8$\\
$M = k[x_5, x_6, x_7, x_8, g]$.\\[0ex]
[If $\lambda = 0$ then $p = x_8^2$ and $f = x_5^2 + x_6^2 - 2x_5x_7 + 2g$].

\item[22.] $\mathfrak{g}_{8.3}$\\
$[x_1, x_j] = x_{j+1}$, $j = 2,\ldots, 7$,\\
$[x_2, x_3] = x_7$, $[x_2, x_4] = x_8$, $[x_2, x_7] = x_8$, $[x_3, x_6] = -x_8$, $[x_4,x_5] = x_8$.\\
$i=2$, $c=5$, $p=x_8$, $F = \langle x_1 - x_2, x_3, x_4, x_5, x_6, x_7, x_8\rangle$, no CP's,\\
$\alpha = 4$, $\mathfrak{h} = C^4 = \langle x_5, x_6, x_7, x_8\rangle$, $Y = k[x_8, f]$,\\
$f = 3x_5^2 x_8 - 2x_7^3 - 6x_5x_8^2 + 6x_6 x_7 x_8 + 6x_8 g$ where\\
$g = (x_1-x_2)x_8 - x_4x_6 + x_3x_7$, $M = k[x_5, x_6, x_7, x_8, g]$.

\item[23.] $\mathfrak{g}_{8.4} \cong Q_8$\\
$[x_1, x_j] = x_{j+1}$, $j = 2,\ldots, 7$,\\
$[x_2, x_7] = x_8$, $[x_3, x_6] = -x_8$, $[x_4,x_5] = x_8$.\\
$i=2$, $c=5$, $p=x_8^2$, $F = \langle x_1 - x_2, x_3, x_4, x_5, x_6, x_7, x_8\rangle$, no CP's,\\
$\alpha = 4$, $\mathfrak{h} = C^4 = \langle x_5, x_6, x_7, x_8\rangle$, $Y = k[x_8, f]$,\\
$f = x_5 ^2 + 2g$ where $g= (x_1-x_2)x_8 + x_3x_7 - x_4x_6$\\
$M = k[x_5, x_6, x_7, x_8, g]$.

\item[24.] $\mathfrak{g}_{8.5} (\lambda \neq 1,2)$\\
$[x_1, x_j] = x_{j+1}$, $j = 2,\ldots, 7$,\\
$[x_2, x_3] = \lambda x_5$, $[x_2, x_4] = \lambda x_6$, $[x_2, x_5] = (\lambda-1) x_7-x_8$,\\
$[x_2, x_6] = (\lambda - 2)x_8$, $[x_3, x_4] = x_7 + x_8$, $[x_3,x_5] = x_8$.\\
$i=2$, $c=5$, $p=1$, $F = \langle x_4, x_5, x_6, x_7, x_8\rangle = C^3 =\mathfrak{h} = CP$, $\alpha = 5$,\\
$Y = k [x_8, f]$, $f = 12x_6x_8^2 + 3(\lambda-1)x_7^4 - 12(\lambda -2) x_4x_8^3 + 6\lambda x_6^2x_8^2$\\
$+12(\lambda - 2)  x_5 x_8^3 + 4(\lambda-2) x_7^3 x_8 - 12(\lambda - 1) x_6 x_7^2x_8 - 12(\lambda -2) x_6x_7 x_8^2$\\
$+12 (\lambda-2) x_5x_7x_8^2 - 6x_7^2 x_8^2$, $M = k[x_4, x_5, x_6, x_7, x_8]$.\\[0ex]
[If $\lambda = 1$ then $p = x_8$ and \\
$f = 6x_6x_8^2 - 6x_5x_8^2 + 6x_4x_8^2 + 3x_6^2 x_8 + 6x_6x_7x_8 - 3x_7^2x_8 - 2x_7^3$].

\item[25.] $\mathfrak{g}_{8.5} (\lambda = 2)$\\
$[x_1, x_j] = x_{j+1}$, $j = 2,\ldots, 7$,\\
$[x_2, x_3] = 2x_5$, $[x_2, x_4] = 2x_6$, $[x_2, x_5] = x_7-x_8$,\\
$[x_3,x_4] = x_7 + x_8$, $[x_3,x_5] = x_8$.\\
$i=2$, $c=5$, $p=x_7^2 -2x_6x_8 - x_8^2$, $F = \langle x_6, x_7, x_8\rangle$,\\
$\alpha = 5$, $\mathfrak{h} = C^3 = \langle x_4, x_5, x_6, x_7, x_8\rangle = CP$,\\
$Y = k [x_8, 2x_6x_8 - x_7^2]$, $M = k[x_4, x_5, x_6, x_7, x_8]$.

\item[26.] $\mathfrak{g}_{8.6} (\lambda \neq 1, 2)$\\
$[x_1, x_j] = x_{j+1}$, $j = 2,\ldots, 7$,\\
$[x_2, x_3] = \lambda x_5$, $[x_2, x_4] = \lambda x_6$, $[x_2, x_5] = (\lambda - 1)x_7$,\\
$[x_2,x_6] = (\lambda - 2)x_8$, $[x_3, x_4] = x_7$, $[x_3,x_5] = x_8$.\\
$i=2$, $c=5$, $p=1$, $F = \langle x_4, x_5, x_6, x_7, x_8\rangle = C^3 = \mathfrak{h} = CP$,\\
$\alpha = 5$,  $Y = k[x_8, f]$\\
$f = 4(\lambda-2) x_4x_8^3 - 2\lambda x_6^2x_8^2 - 4(\lambda - 2) x_5 x_7x_8^2 + 4(\lambda -1) x_6x_7^2 x_8 - (\lambda-1) x_7^4$,\\
$M = k[x_4, x_5, x_6, x_7, x_8]$.\\[0ex]
[If $\lambda = 1$ then $p = x_8^2$, $f = 2x_4x_8 + x_6^2 -2x_5x_7$].

\item[27.] $\mathfrak{g}_{8.6} (\lambda = 2)$\\
$[x_1, x_j] = x_{j+1}$, $j = 2,\ldots, 7$,\\
$[x_2, x_3] = 2x_5$, $[x_2, x_4] = 2x_6$, $[x_2, x_5] = x_7$, $[x_3,x_4] = x_7$, $[x_3,x_5] = x_8$,\\
$i=2$, $c=5$, $p=2x_6x_8 - x_7^2$, $F = \langle x_6, x_7, x_8\rangle$, $\alpha = 5$,\\
$\mathfrak{h} = C^3 = \langle x_4, x_5, x_6, x_7, x_8\rangle = CP$,\\
$Y = k [x_8, 2x_6x_8 - x_7^2]$, $M = k[x_4, x_5, x_6, x_7, x_8]$.

\item[28.] $\mathfrak{g}_{8.7}$\\
$[x_1, x_j] = x_{j+1}$, $j = 2,\ldots, 7$,\\
$[x_2, x_3] = x_6$, $[x_2, x_4] = x_7$, $[x_2, x_5] = -x_7+ x_8$, $[x_2,x_6] = -2x_8$,\\
$[x_3, x_4] = x_7$, $[x_3, x_5] = x_8$,\\
$i=2$, $c=5$, $p=1$, $\alpha = 5$, $F = \langle x_4, x_5, x_6, x_7, x_8\rangle = C^3 = \mathfrak{h} = CP$,\\
$Y = k [x_8, f]$, $f = 8x_4 x_8^3 - 8x_5x_7x_8^2+ 4x_6x_7^2x_8 - x_7^4$,\\
$M = k[x_4, x_5, x_6, x_7, x_8]$.

\item[29.] $\mathfrak{g}_{8.8} (\lambda)$\\
$[x_1, x_j] = x_{j+1}$, $j = 2,\ldots, 7$,\\
$[x_2, x_3] = x_5+ \lambda x_7$, $[x_2, x_4] = x_6 + \lambda x_8$, $[x_2, x_5] = x_7 - x_8$, \\
$[x_2, x_6] = x_8$, $[x_3,x_4] = x_8$,\\
$i=2$, $c=5$, $p=x_8$, $F = \langle x_5, x_6, x_7, x_8\rangle$, $\alpha = 5$,\\
$\mathfrak{h} = C^3 = \langle x_4, x_5, x_6, x_7, x_8\rangle = CP$,\\
$Y = k [x_8, f]$, $f = 6(x_5 + x_6)x_8^2 - 6x_6x_7x_8 - 3x_7^2x_8 + 2x_7^3$,\\
$M = k[x_4, x_5, x_6, x_7, x_8]$.

\item[30.] $\mathfrak{g}_{8.9} (\lambda \neq 1)$ (see Example 9.1)\\
$[x_1, x_j] = x_{j+1}$, $j = 2,\ldots, 7$,\\
$[x_2, x_3] = \lambda x_6+ x_7$, $[x_2, x_4] = \lambda x_7 + x_8$, $[x_2, x_5] = (\lambda - 1)x_8$, \\
$[x_3, x_4] = x_8$,\\
$i=2$, $c=5$, $p=x_8^2$, $F = \langle x_6, x_7, x_8\rangle$, $\alpha = 5$,\\
$\mathfrak{h} = C^3 = \langle x_4, x_5, x_6, x_7, x_8\rangle = CP$,\\
$Y = k [x_8, 2x_6x_8 - x_7^2]$, $M = k[x_4, x_5, x_6, x_7, x_8]$.

\item[31.] $\mathfrak{g}_{8.10} (\lambda \neq 1)$\\
$[x_1, x_j] = x_{j+1}$, $j = 2,\ldots, 7$,\\
$[x_2, x_3] = \lambda x_6$, $[x_2, x_4] = \lambda x_7 $, $[x_2, x_5] = (\lambda - 1)x_8$, $[x_3,x_4] = x_8$.\\
$i=2$, $c=5$, $p=x_8^2$, $F = \langle x_6, x_7, x_8\rangle$, $\alpha = 5$,\\
$\mathfrak{h} = C^3 = \langle x_4, x_5, x_6, x_7, x_8\rangle = CP$,\\
$Y = k [x_8, 2x_6x_8 - x_7^2]$, $M = k[x_4, x_5, x_6, x_7, x_8]$.
\end{itemize}
\newpage

{\bf 8.b {$\boldmath{\mathfrak{g}}$} is not coregular (i.e. {$\boldmath{i(g) \geq  4}$})}
\begin{itemize}
\item[32.] $\mathfrak{g}_{8.9}$ (see Example 9.1)\\
$[x_1, x_j] = x_{j+1}$, $j = 2,\ldots, 7$,\\
$[x_2, x_3] = x_6+ x_7$, $[x_2, x_4] = x_7 + x_8$, $[x_3, x_4] = x_8$, \\
$i=4$, $c=6$, $p=1$, $F = \langle x_2, x_3, x_4, x_5, x_6, x_7, x_8\rangle$, no CP's, $\alpha = 5$,\\
$\mathfrak{h} = C^3 = \langle x_4, x_5, x_6, x_7, x_8\rangle$, $Q(Y) = k(x_8, f_1, f_2, f_3)$,\\
$f_1 = 2x_6x_8 - x_7^2$, $f_2 = 3x_5x_8^2 - 3x_6x_7x_8 + x_7^3$,\\
$f_3 = 10(x_2 - x_3) x_8^2 - 10x_3x_7x_8 + 10x_4x_6x_8 + 10x_4x_7x_8 - 5x_5^2x_8 - 2x_5x_6x_8 + 2x_6^2x_7 - 4x_5x_7^2$,\\
$M = k[x_4, x_5, x_6, x_7, x_8, x_2x_8 - x_3x_8 - x_3x_7]$.

\item[33.] $\mathfrak{g}_{8.10} (\lambda = 1)$\\
$[x_1, x_j] = x_{j+1}$, $j = 2,\ldots, 7$,\\
$[x_2, x_3] = x_6$, $[x_2, x_4] = x_7 $, $[x_3, x_4] = x_8$, \\
$i=4$, $c=6$, $p=1$, $F = \langle x_2, x_3, x_4, x_5, x_6, x_7, x_8\rangle$, no CP's, $\alpha = 5$,\\
$\mathfrak{h} = C^3 = \langle x_4, x_5, x_6, x_7, x_8\rangle$, $Q(Y) = k(x_8, f_1, f_2, f_3)$,\\
$f_1 = 2x_6x_8 - x_7^2$, $f_2 = 2x_2x_8 - 2x_3x_7 + 2x_4x_6 - x_5^2$,\\
$f_3 = 3x_5x_8^2 - 3x_6x_7x_8 + x_7^3$,\\
$M = k[x_4, x_5, x_6, x_7, x_8, x_2x_8 - x_3x_7]$.

\item[34.] $\mathfrak{g}_{8.11} (\lambda \neq 0)$\\
$[x_1, x_j] = x_{j+1}$, $j = 2,\ldots, 7$,\\
$[x_2, x_3] = \lambda x_5 + x_7 + x_8$, $[x_2, x_4] = \lambda x_6 + x_8$, $[x_2, x_5] = \lambda x_7$,\\
$[x_2, x_6] = \lambda x_8$,\\
$i=4$, $c=6$, $p=1$, $F = \langle x_3, x_4, x_5, x_6, x_7, x_8\rangle = C^2 = \mathfrak{h} = CP$, $\alpha = 6$,\\
$Q(Y) = k(x_8, f_1, f_2, f_3)$,\\
$f_1 = 3x_5x_8^2 - 3x_6x_7x_8 + x_7^3$,\\
$f_2 = 4(x_6 - \lambda x_4) x_8^3 + 2\lambda (x_6^2 + 2x_5x_7)x_8^2 - 2x_7^2x_8^2 - 4\lambda x_6 x_7^2x_8 + \lambda x_7^4$,\\
$f_3 = 10(x_6 - \lambda x_3) x_8^4 - 5x_7^2x_8^3 + 2\lambda (5x_4x_7x_8^3 + 5x_5 x_6 x_8^3 - 5x_5x_7^2x_8^2 - 5x_6^2x_7x_8^2 + 5x_6x_7^3x_8 - x_7^5)$,\\
$M = k[x_3, x_4, x_5, x_6, x_7, x_8]$.\\[0ex]
[If $\lambda = 0$ then $Q(Y) = k(x_8, f, g, h)$ with \\
$f = 2x_6x_8 - x_7^2$, $g = 3x_5x_8^2 - 3x_6x_7x_8 + x_7^3$,\\
$h = 10(x_3-x_4) x_8^2 - 2x_5x_6x_8 - 10x_4x_7x_8 - 5x_6^2x_8 + 10x_5x_7x_8 + 4x_5x_7^2 - 2x_6^2x_7$].

\item[35.] $\mathfrak{g}_{8.12}$\\
$[x_1, x_j] = x_{j+1}$, $j = 2,\ldots, 7$,\\
$[x_2, x_3] = x_5 + x_7$, $[x_2, x_4] = x_6 + x_8$, $[x_2, x_5] = x_7$, $[x_2, x_6] = x_8$.\\
$i=4$, $c=6$, $p=1$, $\alpha = 6$, $F = \langle x_3, x_4, x_5, x_6, x_7, x_8\rangle = C^2 = \mathfrak{h} = CP$,\\
$Q(Y) = k(x_8, f_1, f_2, f_3)$, $f_1 = 3x_5x_8^2 - 3x_6x_7x_8 + x_7^3$,\\
$f_2 = 4(x_4 - x_6) x_8^3 - 2 x_6^2x_8^2 + 2x_7^2x_8^2 + 4x_6 x_7^2x_8 - 4x_5x_7x_8^2 - x_7^4$,\\
$f_3 = 5x_3x_8^4 - 5x_5x_6x_8^3 - 5x_4 x_7 x_8^3 + 5x_5x_7^2x_8^2+ 5x_6^2x_7x_8^2 - 5x_6x_7^3x_8 + x_7^5)$,\\
$M = k[x_3, x_4, x_5, x_6, x_7, x_8]$.

\item[36.] $\mathfrak{g}_{8.13}$\\
$[x_1, x_j] = x_{j+1}$, $j = 2,\ldots, 7$,\\
$[x_2, x_3] = x_5 + x_8$, $[x_2, x_4] = x_6$, $[x_2, x_5] = x_7$, $[x_2, x_6] = x_8$.\\
$i=4$, $c=6$, $p=1$, $\alpha = 6$, $F = \langle x_3, x_4, x_5, x_6, x_7, x_8\rangle = C^2 = \mathfrak{h} = CP$,\\
$Q(Y) = k(x_8, f_1, f_2, f_3)$, $f_1 = 3x_5x_8^2 - 3x_6x_7x_8 + x_7^3$,\\
$f_2 = 4x_4x_8^3 - 4x_5x_7x_8^2 - 2x_6^2x_8^2 + 4x_6x_7^2x_8 - x_7^4$,\\
$f_3 = 30(x_3 - x_6) x_8^4 +32x_4^2x_8^3 + 15x_7^2x_8^3 - 30x_3x_5x_8^3 - 30x_4x_7x_8^3 - 34 x_4x_5x_7x_8^2$\\
$+30x_3x_6x_7x_8^2 + 15x_5x_7^2x_8^2 + 30x_5^2x_6x_8^2 - 32x_4x_6^2x_8^2 - 5x_6x_7^3x_8 + 8x_6^4x_8$\\
$-28x_5 x_6^2 x_7 x_8 + 34x_4x_6x_7^2x_8 - 10x_3x_7^3x_8 + 2x_5^2x_7^2x_8 - 2x_6^3x_7^2 + 6x_5x_6x_7^3$\\
$-2x_6^3x_7^2 - 6x_4x_7^4 + x_7^5$.\\
$M = k[x_3, x_4, x_5, x_6, x_7, x_8]$.

\item[37.] $\mathfrak{g}_{8.14} = R_8$\\
$[x_1, x_j] = x_{j+1}$, $j = 2,\ldots, 7$,\\
$[x_2, x_3] = x_5$, $[x_2, x_4] = x_6$, $[x_2, x_5] = x_7$, $[x_2, x_6] = x_8$.\\
$i=4$, $c=6$, $p=1$, $\alpha = 6$, $F = \langle x_3, x_4, x_5, x_6, x_7, x_8\rangle = C^2 = \mathfrak{h} = CP$,\\
$Q(Y) = k(x_8, f_1, f_2, f_3)$, $f_1 = 3x_5x_8^2 - 3x_6x_7x_8 + x_7^3$, \\
$f_2 = 4x_4x_8^3 - 4x_5x_7x_8^2 - 2x_6^2x_8^2 + 4x_6x_7^2x_8 - x_7^4$,\\
$f_3 = 5x_3x_8^4 - 5x_5x_6x_8^3 - 5x_4x_7x_8^3 + 5x_6^2x_7x_8^2 + 5x_5x_7^2x_8^2 - 5x_6x_7^3x_8 + x_7^5$,\\
$M = k[x_3, x_4, x_5, x_6, x_7, x_8]$.

\item[38.] $\mathfrak{g}_{8.15}$\\
$[x_1, x_j] = x_{j+1}$, $j = 2,\ldots, 7$,\\
$[x_2, x_3] = x_6+x_7$, $[x_2, x_4] = x_7+x_8$, $[x_2, x_5] = x_8$.\\
$i=4$, $c=6$, $p=1$, $\alpha = 6$, $F = \langle x_3, x_4, x_5, x_6, x_7, x_8\rangle = C^2 = \mathfrak{h} = CP$,\\
$Q(Y) = k(x_8, f_1, f_2, f_3)$, $f_1 = 2x_6x_8 - x_7^2$,\\
$f_2 = 6(x_4 - x_5)x_8^2 + 3x_6^2x_8 - 6x_5x_7x_8 + 6x_6x_7x_8 - 2x_7^3$,\\
$f_3 = 10x_3x_8^3 - 12x_4x_6x_8^2 + 2x_5x_6x_8^2 - 10x_4x_7x_8^2 + 12x_5x_6x_7x_8$,\\
$-2x_6^2x_7x_8 - 6x_6^3x_8 + 6x_4x_7^2x_8 + 4x_5x_7^2x_8 + 3x_6^2x_7^2 - 6x_5x_7^3$,
$M = k[x_3, x_4, x_5, x_6, x_7, x_8]$.

\item[39.] $\mathfrak{g}_{8.16}$\\
$[x_1, x_j] = x_{j+1}$, $j = 2,\ldots, 7$,\\
$[x_2, x_3] = x_6$, $[x_2, x_4] = x_7$, $[x_2, x_5] = x_8$.\\
$i=4$, $c=6$, $p=1$, $\alpha = 6$, $F = \langle x_3, x_4, x_5, x_6, x_7, x_8\rangle = C^2 = \mathfrak{h} = CP$,\\
$Q(Y) = k(x_8, f_1, f_2, f_3)$, $f_1 = 2x_6x_8 - x_7^2$, 
$f_2 = 2x_4x_8 + x_6^2 - 2x_5x_7$,\\
$f_3 = 5x_3x_8^4 - 5x_4x_7x_8^3 - 5x_5x_6x_8^3 + 5x_5x_7^2x_8^2 + 5x_6^2x_7x_8^2 - 5x_6x_7^3x_8 + x_7^5$,\\
$M = k[x_3, x_4, x_5, x_6, x_7, x_8]$.

\item[40.] $\mathfrak{g}_{8.17}$\\
$[x_1, x_j] = x_{j+1}$, $j = 2,\ldots, 7$,\\
$[x_2, x_3] = x_7$, $[x_2, x_4] = x_8$.\\
$i=4$, $c=6$, $p=1$, $\alpha = 6$, $F = \langle x_3, x_4, x_5, x_6, x_7, x_8\rangle = C^2 = \mathfrak{h} = CP$,\\
$Q(Y) = k(x_8, f_1, f_2, f_3)$, $f_1 = 2x_6x_8 - x_7^2$, 
$f_2 = 3x_5x_8^2 - 3x_6x_7x_8 + x_7^3$,\\
$f_3 = 5x_3x_8^2 + x_5x_6x_8 - 5x_4x_7x_8 - x_6^2x_7 + 2x_5x_7^2$,\\
$M = k[x_3, x_4, x_5, x_6, x_7, x_8]$.

\item[41.] $\mathfrak{g}_{8.18}$\\
$[x_1, x_j] = x_{j+1}$, $j = 2,\ldots, 7$,\\
$[x_2, x_3] = x_8$.\\
$i=4$, $c=6$, $p=x_8$, $F = \langle x_4, x_5, x_6, x_7, x_8\rangle$, $\alpha = 6$,\\
$\mathfrak{h} = C^2 = \langle x_3, x_4, x_5, x_6, x_7, x_8\rangle = CP$,\\
$Q(Y) = k(x_8, f_1, f_2, f_3)$, $f_1 = 2x_6x_8 - x_7^2$, 
$f_2 = 2x_4x_8 + x_6^2 - 2x_5x_7$,\\
$f_3 = 3x_5x_8^2 - 3x_6x_7x_8 + x_7^3$,\\
$M = k[x_3, x_4, x_5, x_6, x_7, x_8]$.

\item[42.] $\mathfrak{g}_{8.19} = L_8$\\
$[x_1, x_j] = x_{j+1}$, $j = 2,\ldots, 7$,\\
$i=6$, $c=7$, $p=1$, $\alpha = 7$, $F = \langle x_2, x_3, x_4, x_5, x_6, x_7, x_8\rangle = \mathfrak{h} = CP$, \\
$Q(Y) = k(x_8, f_1, f_2, f_3)$, $f_1 = 2x_6x_8 - x_7^2$, 
$f_2 = 3x_5x_8^2 - 3x_6x_7x_8 + x_7^3$,\\
$f_3 = 2x_4x_8 + x_6^2 - 2x_5x_7$,\\
$f_4 = 5x_3x_8^2 + x_5x_6x_8 - 5x_4x_7x_8 - x_6^2x_7 + 2x_5x_7^2$,\\
$f_5 = 2x_2 x_8 + 2x_4x_6 - x_5^2 - 2x_3x_7$\\
$M = k[x_2, x_3, x_4, x_5, x_6, x_7, x_8]$.
\end{itemize}

{\bf Corollary 9.3}\\
The Milovanov conjecture holds for all complex filiform Lie algebras of dimension $\leq 8$.  Among these Lie algebras only 9 (namely 5, 8, 18, 20, 21, 22, 23, 32, 33) do not possess a CP.\\
\newpage
{\bf Acknowledgments}\\
We would like to thank Rudolf Rentschler for his inspiring observation and for his helpful comments. We are also very grateful to Manuel Ceballos for his kind offer to compute $\alpha(L_{8,16})$ and $\alpha(L_{8,18})$
with his algorithm.  Furthermore, we thank Alexander Elashvili for suggesting to include the metabelian case.  Finally, we wish to thank Viviane Mebis for her excellent typing of the manuscript.\\

\end{document}